\documentclass[11pt,english]{article}
\usepackage{graphicx} 
\usepackage{amssymb, amsthm, amsmath, fullpage, hyperref}
\usepackage{graphicx, enumerate}
\usepackage{tikz-network}
\usetikzlibrary{positioning, 
                quotes}
\usepackage[utf8]{inputenc}
\usepackage{cleveref}
\usepackage{amsfonts}
\usepackage{amssymb}
\usepackage{amsmath,amsthm,dsfont}
\usepackage[english]{babel}
\usepackage{makeidx}
\usepackage{mathtools}
\usepackage{bbm}
\usepackage{fullpage,amsfonts,amssymb,epsfig,epstopdf,amsmath,titling,url,array,titlesec}
\usepackage{mathtools}
\usepackage{authblk}
\usepackage{tikz}
\usepackage{float}
\usepackage{makeidx}
\usepackage[mathscr]{euscript}
\usepackage{nomencl}
\usepackage{comment}
\theoremstyle{plain}
\newtheorem{thm}{Theorem}[section]
\newtheorem{lem}[thm]{Lemma}
\newtheorem{cor}[thm]{Corollary}

\theoremstyle{definition}
\newtheorem{defn}[thm]{Definition}

\newtheorem{rem}[thm]{Remark}

\newtheorem{claim}[thm]{Claim}
\newtheorem*{notn}{Notation}
\providecommand{\keywords}[1]
{
  \textbf{\textit{Keywords:}} #1
}

\author[1]{Rajat Adak}
\author[2]{L. Sunil Chandran}
\affil[1]{\texttt{rajatadak@iisc.ac.in}} \affil[2]{\texttt{sunil@iisc.ac.in}}
\affil[1,2]{Department of Computer Science and Automation,
Indian Institute of Science, Bangalore, India}
\date{}

\title{Vertex-Based Localization of Erd\H{o}s-Gallai Theorems for Paths and Cycles}

\begin{document}

\maketitle
\begin{abstract}
For a simple graph $G$, let $n$ and $m$ denote the number of vertices and edges in $G$, respectively. The Erd\H{o}s-Gallai theorem for paths states that in a simple $P_k$-free graph, $m \leq \frac{n(k-1)}{2}$, where $P_k$ denotes a path with length $k$ (that is, with $k$ edges). 
In this paper, we generalize this result as follows: For each $v \in V(G)$, let $p(v)$ be the length of the longest path that contains $v$. We show that  
\[m \leq \sum_{v \in V(G)} \frac{p(v)}{2}\]  
The Erd\H{o}s-Gallai theorem for cycles states that in a simple graph $G$ with circumference (that is, the length of the longest cycle) at most $k$, we have  
$m \leq \frac{k(n-1)}{2}$.
We strengthen this result as follows: For each $v \in V(G)$, let $c(v)$ be the length of the longest cycle that contains $v$, or $2$ if $v$ is not part of any cycle. We prove that  
\[m \leq \left( \sum_{v \in V(G)} \frac{c(v)}{2} \right) - \frac{c(u)}{2}\]  
where $c(u)$ denotes the circumference of $G$.  
\newline Furthermore, we characterize the class of extremal graphs that attain equality in these bounds.

\end{abstract}

\section{Introduction}
Typical problems in extremal graph theory aim to maximize or minimize the number of edges in a graph, adhering to certain structural restrictions and characterize the extremal graphs achieving the optimal values. Some of the earliest classical problems in extremal graph theory include: 
\begin{thm}\label{th:Turan}\emph{(Tur\'{a}n \cite{turan1941egy})}
    For a simple graph $G$ with $n$ vertices and clique number at most $r$,
    \begin{equation*}
        |E(G)| \leq \frac{n^2(r-1)}{2r}
    \end{equation*}
    and equality holds if and only if $G$ is a regular Tur\'{a}n Graph with $n$ vertices and $r$ classes.
\end{thm}

\begin{thm}\label{thm:path}\emph{(Erd\H{o}s-Gallai \cite{gallai1959maximal})}
For a simple graph $G$ with $n$ vertices without a path of length $k (\geq 1)$, that is a path with $k$ edges,
\begin{equation*}
    |E(G)| \leq \frac{n(k-1)}{2}
\end{equation*}
and equality holds if and only if all the components of $G$ are complete graphs of order $k$.
\end{thm}

\begin{thm}\label{th:clas}\emph{(Erd\H{o}s-Gallai \cite{gallai1959maximal})}
    For a simple graph $G$ with $n$ vertices and circumference at most $k$,
    \[|E(G)| \leq \frac{k(n-1)}{2}\]
    and equality holds if and only if $G$ is connected and all its blocks are complete graphs of order $k$.
\end{thm}

\subsection{Localization}
A well-known bound on the independence number, $\alpha(G)$, of a graph $G$ is given by: 
\begin{equation}\label{ind} \alpha(G) \geq \frac{n}{\Delta +1} \end{equation} 
where $\Delta$ represents the maximum degree of a vertex in $G$. Independently, Caro \cite{caro1979new} and Wei \cite{wei1981lower} generalized this result by replacing $\Delta$, a global graph parameter, with the degree of individual vertices, which is a \textit{local parameter} for the vertices. The bound they proposed is: 
\begin{equation}\label{indg} \alpha(G) \geq \sum_{v \in V(G)} \frac{1}{d(v)+1} \end{equation} 
where $d(v)$ denotes the degree of a vertex $v \in V(G)$. This leads to the question of whether the aforementioned extremal results can be generalized by considering some localized vertex parameters. For instance, could we replace the global parameter, namely the circumference of the graph, in \Cref{th:clas}, with the length of the largest cycle that each vertex is part of, which is a vertex-based local parameter?
\newline  Localizing the parameters based on edges has already been studied in the literature. But the way they generalized the classical theorems is not as direct as the generalization of \ref{ind} to \ref{indg} established by Caro and Wei.
\vspace{2mm}
\newline For example, Bradač \cite{bradac} proposed a generalization of Turán's \Cref{th:Turan} by associating weights to the edges of the graph. The edge weights were defined as follows:
\begin{equation*}
    k(e) = max\{r \mid e\ \text{occurs in a subgraph of $G$ isomorphic to}\ K_r\}
\end{equation*}
Subsequently, Malec and Tompkins \cite{DBLP:journals/ejc/MalecT23} provided an alternative proof for this generalized version of Turán's theorem.
\begin{thm}\label{thm:GTuran}\emph{(Bradač \cite{bradac} and Malec-Tompkins \cite{DBLP:journals/ejc/MalecT23})} For a simple graph $G$ with $n$ vertices,
\begin{equation*}
    \sum_{e \in E(G)}\frac{k(e)}{k(e) -1} \leq \frac{n^2}{2}
\end{equation*}
and equality holds if and only if $G$ is a regular Tur\'{a}n graph.
\end{thm}
\begin{rem}
    If $G$ is $K_{r+1}$ free, that is $k(e) \leq r$ for all $e \in E(G)$, replacing $k(e)$ by $r$ in \Cref{thm:GTuran} yields the classical Tur\'{a}n bound as stated in \Cref{th:Turan} since $\frac{r}{r-1} \leq \frac{k(e)}{k(e)-1}$ for all $e \in E(G)$.
\end{rem}
 \noindent Malec and Tompkins \cite{DBLP:journals/ejc/MalecT23} coined the term \textit{localization} to describe such weight assignments, as the weight of each edge depends solely on the local structures in which it is involved. Recently, Kirsch and Nir \cite{Kirsch_2024} extended \Cref{thm:GTuran} by localizing an extension of \Cref{th:Turan} originally introduced by Zykov \cite{Zykov1949-qr}. Arag\~{a}o and Souza \cite{aragao2024localised} further advanced the theory by developing a localized version of the Graph Maclaurin Inequalities \cite{Khadzhiivanov1977-xg}.
 \newline Using a similar approach as above, Malec and Tompkins \cite{DBLP:journals/ejc/MalecT23} also developed a \textit{localized} version of \Cref{thm:path}. They defined the edge weights as follows:
\begin{equation*}
    l(e) = max\{r \mid e\ \text{occurs in a subgraph of $G$ isomorphic to}\ P_r\}
\end{equation*}
where $P_r$ is a path with $r$ edges.
\begin{thm}\label{thm:Gpath}\emph{(Malec-Tompkins \cite{DBLP:journals/ejc/MalecT23})} For a simple graph $G$ with $n$ vertices,
\begin{equation*}
    \sum_{e\in E(G)}\frac{1}{l(e)} \leq \frac{n}{2}
\end{equation*}  
and equality holds if and only if each component of $G$ is a clique.
\end{thm}
\begin{rem}
     If $G$ is $P_k$-free, then $l(e) \leq k-1$ for all $e \in E(G)$. Substituting $l(e)$ with $k-1$ in \Cref{thm:Gpath} yields the bound presented in \Cref{thm:path}. Note that for $k=1$, we have $m = 0$.
\end{rem}

\begin{defn}
    Recall that a block in a graph is a maximal connected subgraph containing no cut vertices. A \textit{block graph} is a type of connected graph in which every block is a clique.
\end{defn}

\begin{figure}[H]
    \centering
    \includegraphics[width=6cm]{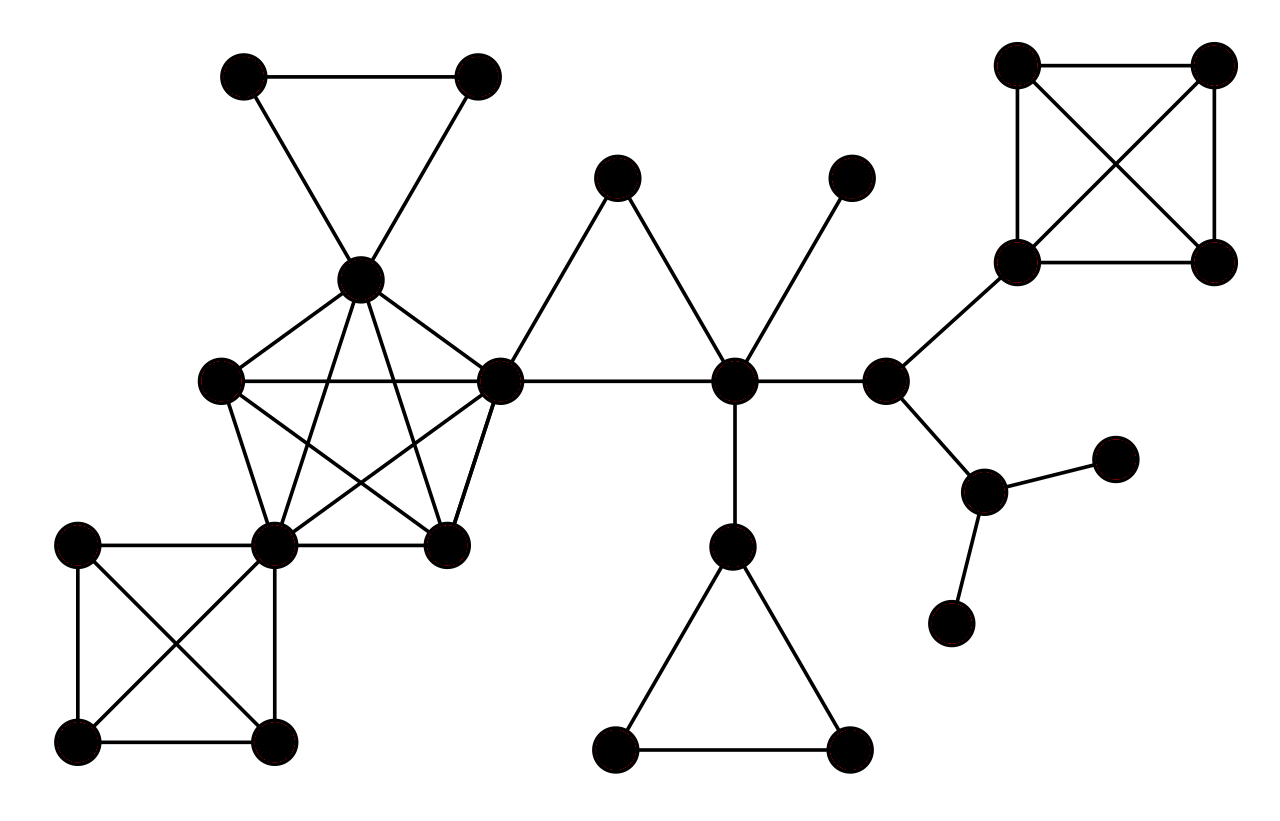}
    \caption{Example of a block graph }
    \label{fig:Path_P}
\end{figure}
\begin{rem}
    Observe that, a tree $T$ can be associated with the block graph $G$ such that the vertex set $V(T)$ corresponds to the blocks of $G$, and two vertices in $T$ are adjacent if and only if their corresponding blocks intersect. Such a tree is referred to as the \textit{block tree} of $G$.

    \noindent Given such a block tree $T$, any vertex in $V(T)$ can be designated as the root, thereby transforming $T$ into a rooted tree with a natural parent-child relationship.
\end{rem}
\begin{defn}
    A block graph $G$ is referred to as a \textit{parent-dominated block graph} if the block tree of $G$ is rooted at a vertex corresponding to a block with largest order, and satisfies the condition that the order of any other block $B$ in $G$ is less than or equal to the order of its parent block.

\end{defn}
\noindent Recently, Zhao and Zhang \cite{https://doi.org/10.1002/jgt.23191} introduced a localized version of \Cref{th:clas}. They assigned a weight to each edge $e \in E(G)$ that is part of some cycle in 
$G$ as follows:
\begin{equation*}
    w(e) = max\{r \mid e\ \text{occurs in a subgraph of $G$ isomorphic to}\ C_r\}
\end{equation*}
For edges that do not belong to any cycle, they set 
$w(e)=2$.
\begin{thm}\emph{(Zhao-Zhang \cite{https://doi.org/10.1002/jgt.23191})}\label{thm:GCycle} For a simple graph $G$ with $n$ vertices,
\begin{equation*}
    \sum_{e\in E(G)} \frac{1}{w(e)} \leq \frac{(n-1)}{2}
\end{equation*}
and equality holds if and only if $G$ is a block graph.
\end{thm}
\begin{rem}
    If the circumference of $G$ is at most $k$, then $w(e) \leq k$ for all $e \in E(G)$, thus substituting $w(e)$ with $k$ in \Cref{thm:GCycle} we retrieve the bound stated in \Cref{th:clas}.
\end{rem}

\section{Our Result: Vertex-based localization}
In Theorems \ref{thm:GTuran}, \ref{thm:Gpath}, and \ref{thm:GCycle}, localization is achieved by assigning weights to the edges of the graph. In this paper, we introduce the concept of \textit{vertex-based localization}, where weights are assigned to the vertices instead of the edges. To the best of our knowledge, this represents the first attempt to achieve localization through assigning weights to the vertices. In this paper, we use vertex-based localization to generalize \Cref{thm:path} and \ref{th:clas}.
\newline While the generalization of \Cref{thm:path} is non-trivial, it is comparatively less challenging than its counterpart, \Cref{th:clas}. The proof of generalization of \Cref{thm:path} follows an inductive approach. In contrast and unfortunately, this technique fails for \Cref{th:clas}, necessitating an entirely new strategy, which significantly increases its complexity and makes the proof more sophisticated. Thus, our main result (in terms of its non-triviality) is \Cref{th:main}, and its proof is presented first, in \cref{Proof of 2.2}. The proof of \Cref{th:main:path} is presented in \cref{Proof of 2.1} and is relatively shorter.
\vspace{2mm}
\newline Following the approach to edge weights introduced by Malec and Tompkins \cite{DBLP:journals/ejc/MalecT23}, we assign a weight to each vertex $v \in V(G)$ as follows:
\begin{equation*}
    p(v) = max\{r \mid v\ \text{occurs in a subgraph of $G$ isomorphic to}\ P_r\} 
\end{equation*}
\begin{thm}\label{th:main:path}
For a simple graph $G$, 
\begin{equation*}
    |E(G)| \leq \sum_{v\in V(G)}\frac{p(v)}{2}
\end{equation*}
Equality holds if and only every connected component of $G$ is a clique.
\end{thm}

\noindent Similar to the edge weights given by Zhao and Zhang \cite{https://doi.org/10.1002/jgt.23191} we assign a weight to each vertex $v \in V(G)$ that is part of a cycle in $G$ as:
\begin{equation*}
    c(v) = max\{r \mid v\ \text{occurs in a subgraph of $G$ isomorphic to}\ C_r\}
\end{equation*}
and for vertices that do not belong to any cycle, we set $c(v) = 2$.
\begin{thm}\label{th:main}
    Given a simple graph, $G$ the number of edges can be upper-bounded as 
\[|E(G)| \leq \left(\sum_{v \in V(G)} \frac{c(v)}{2}\right) - \frac{c(u)}{2}\]
where $c(u)$ is the circumference of $G$. Equality holds if and only if $G$ is a parent-dominated block graph.
\end{thm} 

\subsection{Recovering the Classical Theorems}

To recover the bound and the extremal graph for \Cref{thm:path} from \Cref{th:main:path} we assume $G$ is a $P_k$-free graph with $n$ vertices. Thus, $p(v) \leq (k-1)$, for all $v \in V(G)$. Consequently, we get;
\begin{equation}\label{recover}
    |E(G)| \leq \sum_{v\in V(G)}\frac{p(v)}{2} \leq \sum_{v \in V(G)}\frac{(k-1)}{2} = \frac{n(k-1)}{2}
\end{equation}
This gives the bound as in \Cref{thm:path}. The first inequality of \cref{recover} becomes an equality if and only if all the components of $G$ are cliques, by \Cref{th:main:path}. The second inequality of \cref{recover} becomes an equality if and only if $p(v) = (k-1)$ for all $v \in V(G)$. Thus, all the components of $G$ must be $k$-cliques.
\vspace{2mm}

\noindent Now we show how to recover \Cref{th:clas} from \Cref{th:main}. Since $G$ is a graph on $n$ vertices and with circumference at most $k$, $c(v) \leq k$ for all $v \in V(G)$. Thus, we get;

\begin{equation}\label{eq:clas}
    |E(G)| \leq \left(\sum_{v \in V(G)} \frac{c(v)}{2}\right) - \frac{c(u)}{2} = \left(\sum_{v \in V(G)\setminus \{u\}} \frac{c(v)}{2}\right) \leq \left(\sum_{v \in V(G)\setminus\{u\}} \frac{k}{2}\right) =\frac{k(n-1)}{2}
\end{equation}
 For equality in the first inequality of \Cref{th:main}, $G$ must be a parent-dominated block graph. The second inequality becomes an equality if and only if $c(v) = k$ for all $v \in V(G)\setminus \{u\}.$ Since $c(u)$ is the circumference of $G$, we get $c(u) = k$. Thus the extremal graph for \Cref{th:clas} must be a block graph with blocks of order $k$. (Note that $G$ is a special case of  a parent-dominated block graph).
\vspace{2mm}

\noindent\textbf{Comparison between our results and results of \cite{DBLP:journals/ejc/MalecT23} and \cite{https://doi.org/10.1002/jgt.23191}}
\vspace{2mm}

\noindent Although both \Cref{thm:Gpath} and \Cref{th:main:path} extend \Cref{thm:path}, and both \Cref{thm:GCycle} and \Cref{th:main} strengthen \Cref{th:clas}, the bounds established in our results are fundamentally different. Specifically, \Cref{th:main:path} and \Cref{th:main} provide upper bounds on $|E(G)|$, and in some sense are direct generalizations of the classical Erd\H{o}s-Gallai theorems, by replacing the global parameters by vertex-based local parameters. In contrast, \Cref{thm:Gpath} and \Cref{thm:GCycle} establish lower bounds on the order of the graph using the distribution of edge weights. They are similar in spirit to the generalization of Sperner's Lemma to LYM inequality, from the theory of  set systems (see chapter 3, \cite{bollobás1986combinatorics}). To the best of our knowledge, our bounds and those presented in \cite{DBLP:journals/ejc/MalecT23} and \cite{https://doi.org/10.1002/jgt.23191} do not imply one another.

\subsection{Proof of \Cref{th:main}}\label{Proof of 2.2}
\subsubsection{Some Notations and Lemmas}
\begin{notn}
    For a path $P = v_0v_1\dots v_k$ in graph $G$, we define,
\begin{equation*}
    P(v_i,v_j) = v_iv_{i+1}\dots v_{j-1}v_j \ \text{for } 0\leq i \leq j \leq k 
\end{equation*}
\begin{equation*}
    P(v_i,v_j) = v_iv_{i-1}\dots v_{j+1}v_j \ \text{for } 0\leq j \leq i \leq k
\end{equation*}
\end{notn}
\noindent A $v_0$-path in a graph $G$ is a path that starts at the vertex $v_0 \in V(G)$. Consider a $v_0$-path $P$ in $G$ that ends at a vertex $v_k$. In this paper, unless otherwise stated, we always consider a $v_0$-path $P$, to start from $v_0$ and move towards the other end vertex $v_k$, which we refer to as the \textit{terminal vertex} of $P$. This assumption about the direction will help in defining the following functions.

Define the distance between two vertices $x, y \in V(P)$, as the number of edges between them on path $P$. Denote the distance between $x$ and $y$ on path $P$ as $dist_P(x,y)$. Note that, $dist_P(x,x) = 0$. 

Let $i \in \mathbb{N}$ and $s \in V(P)$ such that $dist_P(v_0, s) \geq i$. Define $Pred_x(i)$ as the vertex $z \in V(P)$ such that $dist_P(z,s) = i$ and $dist_P(v_0,z) \leq dist_P(v_0,s)$, that is, $z$ is a \textit{predecessor} of $s$ on $P$. Similarly, let $j \in \mathbb{N}$ and $t \in V(P)$ such that $dist_P(t,v_k) \geq j$, that is, $z$ is a \textit{successor} of $t$ on $P$. Let $Succ_x(j)$ represent the vertex $z\in V(P)$ such that $dist(t,z) = j$ and $dist_P(t,v_k) \geq dist_P(z,v_k)$.

Let $P= v_0v_1v_2\dots v_k$ be a longest $v_0$-path in $G$. It is easy to see that $N(v_k) \subseteq \{v_0, v_1, \dots, v_{k-1}\}$; otherwise the path $P$ can be extended to get a longer $v_0$-path in $G$. If $v_k$ is adjacent to $v_j$ for some $j \in \{ 0 , 1 \dots k-2\}$, then we can construct an alternate path $P'$ such that $V(P) = V(P')$ and $P'= v_0v_1 \dots v_jv_kv_{k-1} \dots v_{j+2}v_{j+1}$, by removing the edge $v_jv_{j+1}$ and adding the edge $v_jv_k$ to $P$.

\begin{defn}
    Let $P$ be a longest $v_0$-path in $G$, where $P = v_0v_1v_2\dots v_k$ and $v_k$ is adjacent to $v_j$, $0\leq j \leq k-2$. Then the path $P' = v_0v_1 \dots v_jv_kv_{k-1} \dots v_{j+2}v_{j+1}$ is called a \textit{simple transform} of $P.$
\begin{figure}[H]
    \centering
    \includegraphics[width=12cm]{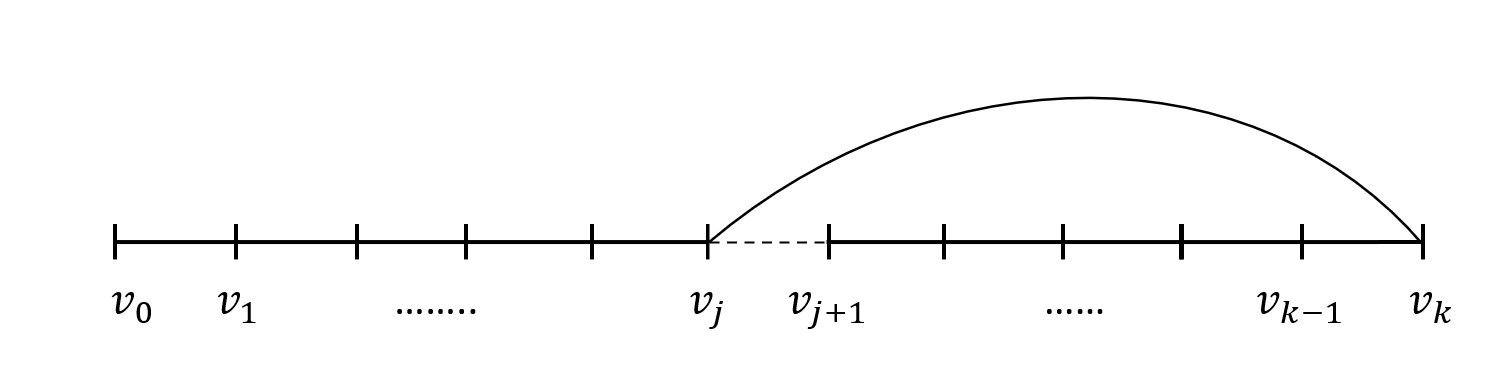}
    \caption{Path $P'$ is a \textit{simple transform} of $P$}
\end{figure}    
\end{defn}
\begin{rem}\label{st_inv}
    If $P'$ is a simple transform of $P$, then $P$ is also a simple transform of $P'$.
\end{rem}
\begin{defn}
   Let $P$ denote a longest $v_0$-path in $G$ with terminal vertex $v_k$. Since, each vertex $x \in N(v_k)$ is on $P$, there exists a simple transform of $P$ ending at $Succ_x(1)$, for each $x \in N(v_k) \setminus \{Pred_x(1)\}$. Let $P'$ be one such simple transform. We can create new paths which are simple transforms of $P'$. Thus by applying a sequence of simple transform operations to $P$, we can obtain several paths. Let $P''$ be one such path. Clearly $P''$ is also a longest $v_0$-path and $V(P'') = V(P)$. We refer to $P''$ as a \textit{transform} of $P$. Let $\mathcal{T}$ represent the set of all transforms of $P$. Note that $\mathcal{T}$ is a function of $(G,P,v_0)$.
\end{defn}
\begin{rem}
    If $P'$ is a transform of $P$, then $P$ is also a transform of $P'$.
\end{rem}
\begin{notn} In a graph $G$ with $P$ as a longest $v_0$-path and $\mathcal{T} = \mathcal{T}(G,P,v_0)$, we define the following notations and symbols, to be used throughout the paper.
\begin{enumerate}
    \item  Let $L(G, P, v_0)$ be the set of terminal vertices of the paths in $\mathcal{T}$; that is: $$L(G,P,v_0) = \{v \in V(P) \mid v \text{ is a terminal vertex of some } P' \in \mathcal{T}\}$$
    \item  Define the set $L^0(G,P,v_0)$ as follows:  
$$L^0(G,P,v_0) = \{v \in V(G) \setminus V(P) \mid \exists\ v' \in L(G,P,v_0) \text{ such that } N(v) = N(v')\}.$$  
It is important to note that the set $L^0(G,P,v_0)$ may be empty. \newline Moreover, $L(G,P,v_0)\cap L^0(G,P,v_0) = \emptyset$. Next we define the set $L^*(G,P,v_0)$ as:  
$$ L^*(G,P,v_0) = L(G,P,v_0) \cup L^0(G,P,v_0).$$

    \item By definition, for $x_1 \in L^0(G,P,v_0)$, there exists $x_2 \in L(G,P,v_0)$ such that $N(x_2) = N(x_1)$. Note that we can get many vertices in $L(G, P,v_0)$ with the same neighborhood as $x_1$. Fix one such vertex for $x_1$, say $x_2$, and refer to $x_2$ as the \textit{twin} of $x_1$, denoted by $tw(x_1) = x_2$. Thus for $x_1 \in L^0(G,P,v_0)$, if there exists a twin of $x_1$, it is unique.
    \item For $v \in L(G,P,v_0)$, let $\mathcal{T}_v$ denote the set of all the transforms in $\mathcal{T}$ with $v$ as the terminal vertex.
     \item For $v \in L^0(G,P,v_0)$, let $x = tw(v)$. We modify all the paths in $\mathcal{T}_x$ by replacing $x$ with $v$ in every path, as illustrated in \Cref{corresponding_path}.
\begin{figure}[H]
    \centering
    \includegraphics[width=12cm]{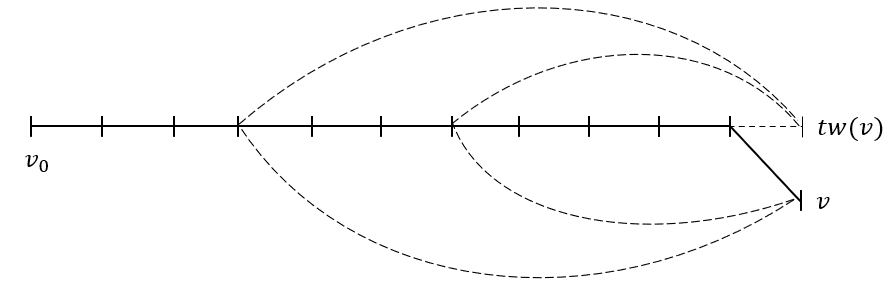}
    \caption{A path from $\mathcal{T}_x$ is modified by replacing $x = tw(v)$ with $v$}
    \label{corresponding_path}
\end{figure}

   We refer to the set of these modified paths as $\mathcal{T}_v$ for some $v \in L^0(G,P,v_0)$. Note that the path in $\mathcal{T}_v$ and $\mathcal{T}_x$ are in one-to-one correspondence.
    \item  For $v \in L^*(G,P,v_0)$, let $S_v$ be the set of neighbors of $v$ outside the set $L(G,P,v_0)$, that is, $S_v = N(v) \cap (L(G,P,v_0))^c$.
         \item Let $w(P_v) = Pred_v(c(v)-1)$ on $P_v \in \mathcal{T}_v$. We call this the \textit{pivot} vertex of the path.  
         \item Define $Front(P_v) = P_v(v_0,Pred_{w(P_v)}(1))$, that is the part of $P_v$ before the pivot. Define $Back(P_v) = P_v(Succ_{w(P_v)}(1),v)$, that is the part of $P_v$ after the pivot. Observe that $v$ is not adjacent to any vertex in $Front(P_v)$, since otherwise we will get a cycle longer than $c(v)$ containing the vertex $v$. 
    \item Define $Front^*(P_v) = P_v(v_0,w(P_v))$ and $Back^*(P_v) = P_v(w(P_v),v)$.
\end{enumerate}
    \end{notn}

\subsubsection{Some Useful Properties}
\begin{lem}\label{same_c(v)_twins}
    For all $v \in L^0(G,P,v_0)$, $c(v) = c(tw(v))$.
    \begin{proof}
        If $c(v) = 2$, vertex $v$ is not part of any cycle in $G$. Since $N(v) = N(tw(v))$, it follows that $|N(v)| =|N(tw(v))|= 1$; otherwise there will be a $4$-cycle containing $v$ and $tw(v)$. Thus $c(tw(v))=2$. Similarly, if we assume $c(tw(v)) = 2$ we get $c(v) =2$.
        
         Now assume $c(v)>2$ and $ c(tw(v)) >2$, that is, both $v$ and $tw(v)$ are part of at least one cycle in $G$. Let $C$ be a cycle of length $c(v)$ that contains the vertex $v$. If the cycle $C$ also contains $tw(v)$, then it follows that $c(v) \leq c(tw(v))$. Otherwise, since $N(v) = N(tw(v))$, we can construct a new cycle $C'$ by replacing $v$ with $tw(v)$ in $C$. This construction yields a cycle of length $c(v)$ containing $tw(v)$, which implies $c(v) \leq c(tw(v))$.

\noindent Similarly, starting with a cycle of length $c(tw(v))$ containing the vertex $tw(v)$, we can deduce that $c(v) \geq c(tw(v))$. These inequalities together establish that $c(v) = c(tw(v))$.
    \end{proof}
\end{lem}
\begin{lem}\label{No_early_nbr}
    Let $c(x) = min\{c(v)\mid v \in L^*(G,P,v_0)\}$. A path $P_v$ is called a good path if the terminal vertex $v$ is not adjacent to any vertex in $P_v(v_0,Pred_vc(x)
    )$; otherwise $P_v$ is a bad path. All the paths in $\mathcal{T}_v$ for all $v \in L^*(G,P,v_0)$ are good.
    \begin{proof}
        Using \Cref{same_c(v)_twins}, without loss of generality we can assume $x\in L(G,P,v_0)$. Let $P_{x} \in \mathcal{T}_{x}$; clearly $x$ can not have a neighbor in $P_x(v_0,Pred_xc(x)) = Front(P_x)$, otherwise there will be cycle of length greater than $c(x)$ containing vertex $x$. Suppose for $v \in L^*(G,P,v_0)$, there exists a bad path $P_v \in \mathcal{T}_v$. If $v\in L^0(G,P,v_0)$, we can work with $tw(v)$ instead of $v$ and the corresponding path $P_{tw(v)}$; thus assume $v \in L(G,P,v_0)$. 
        
        We can get $P_v$ from $P_{x}$, applying a sequence of simple transforms. Assume $P_v$ to be the first bad path in the sequence. Observe that, in all the paths in the sequence till $P_v$ only the position of the final $c(x)-1$ vertices are getting permuted. Thus, the final $c(x) -1$ vertices in all the paths in the sequence starting from $P_x$ to $P_v$ are the same. That is, $V(P_v(z_v,v)) = V(P_v(z_x,x)) = V(Back(P_x))$, where $z_y = Pred_y(c(x)-2)$ on path $P_y \in \mathcal{T}_y$. Moreover, the remaining vertices have invariant positions on the paths. To be more precise, let $i = |V(P)| - c(x) +1$; all the paths in the sequence from $P_x$ to $P_v$ are identical till the $i^{th}$ position starting from $v_0$ in the first position.
        
        Let $v$ be adjacent to $t \in P_v(v_0,Pred_vc(x)
    )$. Thus, the cycle formed by adding the edge $(v,t)$ to path $P_v(t,v)$, is of length more than $c(x)$ and contains the vertex $x$, which is a contradiction.
    \end{proof}
    \end{lem}
     \begin{rem}\label{till i}
     Let $i = |V(P)| - c(x) +1$ where $c(x) = min\{c(v) :\ v\in L^*(G,P,v_0)\}$. As a consequence of \cref{No_early_nbr}, we get that, for any $v_1, v_2 \in L^*(G,P,v_0)$; the paths $P_{v_1} \in \mathcal{T}_{v_1}$ and $P_{v_2} \in \mathcal{T}_{v_2}$ are identical till the $i^{th}$ position starting from $v_0$ in the first position. This is because when a simple transform is formed for $P_v$, the resulting path differs from $P_v$ only after a neighbor of $v$. Thus for $v \in L^*(G,P,v_0)$, the pivot vertex, $w(P_v)$ remains invariant irrespective of $P_v$ chosen from $\mathcal{T}_v$. That is, it does not depend on the path chosen, but only on the (weight of the) terminal vertex. Thus we denote $w(P_v)$ simply by $w_v$.
\end{rem}
    \begin{cor}\label{L in back}
         For all $v \in L^*(G,P,v_0)$ and $P_v \in \mathcal{T}_v$, $L^*(G,P,v_0) \cap V(Front^*(P_v)) = \emptyset$
        \begin{proof} Suppose there exists $v \in L^*(G,P,v_0)$ such that; \[L^*(G,P,v_0) \cap V(Front^*(P_v)) \neq \emptyset \] \[\implies L(G,P,v_0) \cap V(Front^*(P_v)) \neq \emptyset\] Let $t \in L(G,P,v_0) \cap (V(Front^*(P_v))$. Let $P_y \in \mathcal{T}_y$ be such that one can obtain $P_t \in \mathcal{T}_t$ as a simple transform of $P_y$. Thus $y$ is adjacent to $Pred_t(1)$ in $P_y$. But from \cref{No_early_nbr}, $y$ can not have a neighbor in $P_v(v_0,Pred_yc(x))$. Thus we arrive at a contradiction.

        \end{proof}
    \end{cor}
    \begin{rem}
                If fact from \cref{No_early_nbr} we get a stronger statement to be true; namely, $$L^*(G,P,v_0) \cap V(P_v(v_0,Pred_v(c(x)-1)) = \emptyset$$Since, $c(x) = min\{c(v)\ |\ v \in L^*(G,P,v_0)\}$, we get $V(Front^*(P_v)) \subseteq  V(P_v(v_0,Pred_v(c(x)-1))$. Thus, $L^*(G,P,v_0) \cap  V(Front^*(P_v)) = \emptyset$.
            \end{rem}
   
    \begin{lem}\label{|L| = |L_v|}

        From \cref{L in back}, it is clear that for all $v \in L^*(G,P,v_0)$ and $P_v \in \mathcal{T}_v$, all the $L^*(G,P,v_0)$ vertices on $P_v$ are in $Back(P_v)$. There are exactly $|L(G,P,v_0)|$ such vertices in $Back(P_v)$.
        \begin{proof} We will consider the following two cases;

        \begin{enumerate}
            \item If $v \in L(G, P, v_0)$, $P_v$ contains all the vertices from $L(G,P,v_0)$ and none from $L^0(G,P,v_0)$. Thus we get;
            \begin{equation}\label{lv=l}
                L^*(G,P,v_0) \cap V(Back(P_v)) = L(G, P, v_0)
            \end{equation}
            \item Whereas if $v \in L^0(G,P,v_0)$, then $P_v$ contains all the vertices from $L(G,P,v_0)$, except $tw(v)$ instead it $v \in L^0(G,P,v_0)$. Thus we get;
            \begin{equation}\label{Lv = l}
                L^*(G,P,v_0)\cap V(Back(P_v)) = L(G,P,v_0)\setminus \{tw(v)\} \cup \{v\}
            \end{equation}
            
        \end{enumerate}  
        Thus from \cref{lv=l,Lv = l}, we get that the number of $L^*(G,P,v_0)$ elements on any $P_v \in \mathcal{T}_v$ for all $v \in L^*(G,P,v_0)$ is $|L(G,P,v_0)|$.
         \end{proof}
    \end{lem}

\subsubsection{Proof of Inequality}\label{Proof of inequality}
Now we are all set to prove the inequality part of \Cref{th:main}.
\begin{proof}
   Let $G$ be a simple graph and $c(u) = max\{c(v) \mid v \in V(G)\}$. Without loss of generality assume $G$ does not have any isolated vertices, as these vertices do not contribute towards the number of edges. Let $P$ be a longest $u$-path in $G$. We will count the number of edges in $G$ having at least one endpoint in $L^*(G, P,u)$. For this, we will find the contribution of each vertex in $L^*(G,  P,u)$ towards this count. 
    
    \noindent Let $P_v \in \mathcal{T}_v$. Vertex $v$ has exactly $|S_v|$ neighbors outside $L^*(G,P,u)$.  Therefore the number of edges having at least one endpoint in $L^*(G,P,u)$ is given by:
    \begin{equation}\label{eq:1}\sum_{v \in L^*(G,P,u)}  \left(|S_v| + \frac{d(v) - |S_v|}{2}\right) = \sum_{v \in L^*(G,P,u)} \frac{d(v) + |S_v|}{2}    \end{equation} 
    
     \noindent Recall that all the neighbors of $v$ are on $Back^*(P_v)$. Since $|V(Back^*(P_v))| = c(v)$, from \cref{|L| = |L_v|} we get,
    \begin{equation}\label{eq:2}
        |S_v| \leq |V(Back^*(P_v))| - |L^*(G,P,u)\cap Back(P_v)| = c(v) -|L(G,P,u)|
    \end{equation}
   Since for every neighbor $x$ of $v$, the immediate successor of $x$ on $P_v$ belongs to $L^*(G,P,u)$, in the view of \cref{|L| = |L_v|}, we have;
    \begin{equation}\label{eq:3}
        d(v) \leq |L^*(G,P,u) \cap V(Back(P_v))| = |L(G,P,u)| 
    \end{equation}
    Using \cref{eq:2,eq:3} we get the following inequalities;
    \begin{equation}\label{eq:4}
        \sum_{v \in L^*(G,P,u)} \frac{d(v) + |S_v|}{2} \leq \sum_{v \in L^*(G,P,u)} \frac{|L(G,P,u)| + c(v) - |L(G,P,u)|}{2} \leq \sum_{v \in L^*(G,P,u)} \frac{c(v)}{2}
    \end{equation}
    Thus from \Cref{eq:4}, the number of edges with at least one endpoint in $L^*(G,P,u)$ is at most;
    \[\sum_{v \in L^*(G,P,u)} \frac{c(v)}{2}\]
    We proceed as per the following algorithm;
    \begin{enumerate}
        \item Set $G_0 = G$ and $L_0 = L^*(G_0,P,u)$ \item Set $G_1'= G_0 \setminus L_0$, remove all the isolated vertices from $G_1'$ to get $G_1$. 
        \item Take $i =1$ and set $x = u$.
        \item While ($V(G_i) \neq \emptyset$)
        \begin{enumerate}
            \item For each  $v \in V(G_i) $, let $ c_i(v)$ represent the weight of vertex $v$ in the context of $G_i$.

            \item If ($x \notin V(G_i))$
            \begin{enumerate}
                \item Update $x$ such that $c(x) = max\{c(v): v \in V(G_i)\}$.
            \end{enumerate}
             \item Let $P_i$ be the longest $x$-path in $G_i$ \item Set $L_i = L^*(G_i, P_i,x)$.
             \item Set $G_{i+1}' = G_i \setminus L_i$
             \item Remove isolated vertices from $G_{i+1}'$ to get $G_{i+1}$.
            \item Set $i = i+1$
        \end{enumerate}
    \end{enumerate}
    $G$ is a finite graph. Since $|L_i| > 0$ for all $i$, some vertices are being removed from the graph at every step. If $v \in V(G_i)$, then $v_i$ is not isolated in $G_i$. Thus, a non-trivial path $P_i$ can be defined and $|L_i| \neq 0$. Thus, $V(G_i)$ becomes empty in a finite number of steps. Let $t = \max\{i : V(G_i) \neq \emptyset\}$.
    \begin{rem}\label{u notin L}
      Let $j = max\{i : u \in V(G_i)\} \implies u \notin L_k$ where $k \leq j-1$. Path $P_j$ starts at the vertex $u$; thus from the definition of $L_j$, $u \notin L_j$. Since $u \notin V(G_{j+1})$ this implies $u \notin L_k$ for all $k \geq j+1$. Thus we get that $u \notin L_i$ for all $i \in [t]$.
    \end{rem}
    \begin{notn}
        Denote the set of edges having at least one endpoint in $L_i$ as $E^*(L_i)$.
    \end{notn}
    \noindent Now it is easy to see that;
    \begin{equation}\label{eq:(8)}
        E(G) = \bigcup_{i=0}^tE^*(L_i)
    \end{equation}
    Also, we know that;
    \begin{equation}\label{eq:(9)}
        |E^*(L_i)| \leq \sum_{v\in L_i}\frac{c_i(v)}{2}
    \end{equation}
    From \cref{u notin L}, recall that $u \notin L_i$ for any $i \in [t]$. Since for all $v \in V(G_i)$, $c_i(v) \leq c(v)$, using (\ref{eq:(8)}) and (\ref{eq:(9)}) we get;
    \begin{equation}\label{main inequation}
        |E(G)| = \sum_{i=0}^t|E^*(L_i)| \leq \sum_{i=0}^t\left(\sum_{v\in L_i}\frac{c_i(v)}{2}\right)\leq \left(\sum_{v \in V(G)} \frac{c(v)}{2}\right) - \frac{c(u)}{2}
    \end{equation}
\end{proof}
\subsubsection{Characterizing the extremal graphs}
The next natural question is to characterize the class that is extremal for the inequality in \Cref{th:main}
\vspace{2mm}
\newline\textbf{$\bullet$ If part: Parent-dominated block graphs are extremal}
\vspace{2mm}

\noindent We begin by assuming that $G$ is a parent-dominated block graph and show that equality holds in \Cref{th:main} for graph $G$. Let $\{B_1, B_2, \dots, B_k\}$ denote the set of blocks in $G$, where $|V(B_i)| = b_i$ for all $i \in [k]$. Without loss of generality, assume that $B_k$ is the root block and thus a block with the maximum order in $G$. It is easy to see that for a vertex $v$, $c(v) = max\{|B_i| \mid v \in B_i\}$. In particular vertices in $B_k$ have the maximum weight. Let $u \in V(B_k)$ be a vertex with the maximum weight in $G$. For each $i \in [k-1]$, let $w_i$ denote the cut-vertex that connects block $B_i$ with its parent block. Furthermore, assume $w_k = u$. First note that,
\begin{equation}\label{eq}
    |E(G)| = \sum_{i = 1}^k |E(B_i)|
\end{equation}
We proceed by counting the number of edges in $B_i$ while keeping track of the weights of the vertices within $B_i$, excluding the vertex $w_i$. For this purpose, we first note that;

\begin{equation}\label{eq:8}
    \sum_{v \in V(B_i)\setminus\{w_i\}}\frac{c(v)}{2} = \sum_{v \in V(B_i)\setminus\{w_i\}}\frac{b_i}{2} = \frac{(b_i-1)b_i}{2} = |E(B_i)|
\end{equation}

\noindent Now, from \cref{eq,eq:8}, we obtain:
\[\left(\sum_{v \in V(G)} \frac{c(v)}{2}\right) - \frac{c(u)}{2} = \sum_{i = 1}^k \left(\sum_{v \in V(B_i)\setminus\{w_i\}}\frac{c(v)}{2}\right) = \sum_{i =1}^k|E(B_i)| = |E(G)|\]
Thus, when 
$G$ is a parent-dominated block graph, equality holds in \Cref{th:main}.
\vspace{2mm}
\newline\textbf{$\bullet$ Only if part: Equality in \Cref{th:main} implies $G$ is Parent-dominated block graph}
\vspace{2mm}

\noindent Now for the other direction, assume that equality holds in the bound of \Cref{th:main}. We will show that $G$ is a parent-dominated block graph. 
\newline Our approach is to first carefully establish the conditions required to hold in the proof of the inequality of \Cref{th:main}. We note that three conditions are required; which are highlighted in bold face below.
\newline For equality in \Cref{th:main}, we must have equality in (\ref{eq:2}) and (\ref{eq:3}) in the proof of the inequality (refer \cref{Proof of inequality}). That is for all $v \in L^*(G,P,u)$ where $u$ is a maximum weight vertex in $G$ and $P$ is a longest $u$-path, we should have;
\boldmath
\begin{equation}\label{c = l +s}
    c(v) = |L(G,P,u)| +|S_v|
\end{equation}
\begin{equation}\label{eq:9}
    d(v) = |L(G,P,u)|
\end{equation} 
\unboldmath

\noindent For equality we must have $\bigcup_{i=0}^tL_i = V(G) \setminus \{u\}$ and the removal of $L_i$ should not decrease the weight of any vertex in $V(G_{i+1}) \setminus \{u\}$. This follows from the tightness of \cref{main inequation}. Therefore, for all $i \in [t]$ and $\forall$ $v \in V(G_i)\setminus \{u\}$;
\boldmath
\begin{equation}\label{c_i(v) = c(v)}
    c_i(v) = c(v) 
\end{equation} 
\unboldmath
\vspace{2mm}
\newline\textbf{$\bullet$ Preprocessing steps: $G$ is connected and $\delta(G) \geq 2$ assumption}

\begin{lem}\label{connected}
    $G$ is connected.
    \begin{proof}
        Suppose $G$ is not connected. Let $G_1$ and $G_2$ be two distinct connected components. Let $v_1\in V(G_1)$ and $v_2\in V(G_2)$. Construct a new graph $G'$ by adding the edge $(v_1,v_2)$ in $G$. Clearly, the weights of the vertices of $G'$ are the same as those of $G$ since no new cycles are formed. Therefore, in \Cref{th:main} the same upper bound holds for $E(G)$ and $E(G')$. But $|E(G')| > |E(G)|$. Thus we can not have equality for $G$ in the bound of \Cref{th:main}.
    \end{proof}
\end{lem}
\begin{rem}\label{|L| >1}
    Suppose there exists a vertex $x \in V(G)$ such that $d(x) =1$ and thus, $c(x) = 2$. Let $H = G \setminus \{x\}$. Clearly $|E(G)| = |E(H)| +1$. Since there exists at least one non-isolated vertex in $G$ other than $x$, without loss of generality we can assume $x \neq u$. Thus we get $$\left(\sum_{v\in V(G)}\frac{c(v)}{2}\right) - \frac{c(u)}{2} = \left(\sum_{v\in V(H)}\frac{c(v)}{2}\right) - \frac{c(u)}{2} +1$$
    So, in \Cref{th:main}, equality occurs for $G$ if and only if equality occurs for $H$. It is enough to show that 
$H$ is a parent-dominated block graph because this implies that $G$ must also be a parent-dominated block graph.
    \newline Thus we can assume that $\delta(G) \geq 2$. This allows us to assume $|L(G,P,u)|\geq 2$. Since otherwise we get $L(G,P,u) = \{v\}$ and $d(v) =1$, which is a contradiction to $\delta(G) \geq 2$.

\end{rem}
\noindent\textbf{$\bullet$ Holes and Hole-free paths}
\begin{notn}
    Let $v \in L^*(G,P,u)$. $N^+_{P_v}(v)$ be the set of immediate successors of the neighbors of $v$ on path $P_v \in \mathcal{T}$. That is, \[N^+_{P_v}(v) = \{Succ_x(1) \text{ on } P_v :\ x \in N(v)\}\]
\end{notn}
\begin{rem}\label{d(v) = |L|}
    Let $v \in L^*(G,P,u)$. For every neighbor $x$ of $v$, we get a corresponding element of $L^*(G,P,u)$, namely $Succ_x(1)$ on path $P_v \in \mathcal{T}_v$. Suppose there exists an element $y \in L^*(G,P,u)\cap V(Back^*(P_v))$ such that $Pred_y(1)$ on $P_v$ is not adjacent to $v$, then $d(v) < |L^*(G,P,u)\cap V(Back^*(P_v))|$. Recalling \cref{|L| = |L_v|}, we get $d(v) < |L^*(G,P,u)\cap V(Back^*(P_v))| = |L^*(G,P,u)\cap V(Back(P_v))|=|L(G,P,u)|$, thus contradicting \cref{eq:9}. Therefore we get; $$L^*(G,P,u)\cap V(Back(P_v)) = N^+_{P_v}(v)$$
\end{rem}
\begin{defn}
    Let $v \in L^*(G,P,u)$ and $P_v \in \mathcal{T}_v$, a vertex $x$ is called a \textit{hole} if $x \in Back^*(P_v)$ and $x \notin N(v)$, and with respect to path $P_v$, $Pred_x(1) \notin N(v)$. That is, $x \in V(Back^*(P_v))$ is a hole if $x \notin N(v) \cup N^+_{P_v}(v)$.
\end{defn}
\begin{lem}\label{nohole}
    For all $v \in L^*(G,P,u)$, there are no holes in $P_v$.
    \begin{proof}
        From \cref{c = l +s}, we have that $|S_v| + |L(G,P,u)| = c(v) =|Back^*(P_v)|$. Since $S_v \cap L^*(G,P,u) = \emptyset$, a vertex in $Back^*(P_v)$ is either in $S_v$ or $L^*(G,P,u)$. Let $x \in V(Back^*(P_v))$ be an hole in $P_v$, that is, \[x \notin N(v) \cup N^+_{P_v}(v )\] Thus, $x$ is not adjacent to $v$ and hence can not be an element of $S_v$. So, $x$ must be an element of $L^*(G,P,u)$. Thus, from \cref{d(v) = |L|} we get; \[x \in L^*(G,P,u) \cap V(Back^*(P_v)) = N^+_{P_v}(v)\] We arrive at a contradiction since we assumed $ x \notin N^+_{P_v}(v)$.
    \end{proof}
\end{lem}
\begin{rem}\label{wv nbr}
    Let $v \in L^*(G,P,u)$, and $P_v \in \mathcal{T}_v$. From \cref{nohole}, we get that there are no holes in $P_v$; and therefore $w_v$ is itself not a hole. Since $w_v \notin N^+_{P_v}(v)$, $w_v$ should be a neighbor of $v$.
    \vspace{4mm}
    
    \noindent \textbf{$\bullet$ Concepts that coincide for all the paths in $\bigcup_{v\in L^*(G,P,u)} \mathcal{T}_v$ for extremal graphs}
    \begin{lem}\label{same c(v)}
        $c(v)$ is same for all $v \in L^*(G,P,u)$.
        \begin{proof}
           From \Cref{No_early_nbr}, for any $v \in L^*(G,P,u)$, vertex $v$ can not have a neighbor in $P_v(u,Pred_vc(x))$, where $P_v \in \mathcal{T}_v$ and $c(x) = min\{c(v)\ |\ v \in L^*(G,P,u)\}$. Thus $c(v) \geq c(x)$. \newline Suppose $c(v) > c(x)$. On path $P_v$, $w_v  = Pred_v(c(v)-1) \in V(P_v(u,Pred_vc(x)))$. From \cref{wv nbr} we have $w_v \in N(v)$. Thus we arrive at a contradiction. Therefore, $c(v) = c(x)$ for all $v \in L^*(G,P,u)$.
        \end{proof}
    \end{lem}
    \begin{rem}
        From \cref{till i} and \cref{same c(v)} we get that for any $v_1,v_2 \in L^*(G,P,u)$, $w_{v_1} = w_{v_2}$; that is, the pivots coincide for all paths in $\bigcup_{v\in L^*(G,P,u)} \mathcal{T}_v$. As the pivot vertex is independent of the terminal vertex of the path, we can denote it as $w$, dropping the subscript. Restating \cref{till i} for the extremal graph $G$, if $v_1, v_2 \in L(G,P,u)$, the paths $P_{v_1} \in \mathcal{T}_{v_1}$ and $P_{v_2} \in \mathcal{T}_{v_2}$ are identical from $u$ to $w$. Moreover, $V(Back(P_{v_1}))= V(Back(P_{v_2}))$. If $v_3 \in L^0(G,P,u),$ the path $P_{v_3} \in \mathcal{T}_{v_3}$ and $P_{v_1}$ are identical from $u$ to $w$ and $V(Back(P_{v_1})) = V(Back(P_{v_3}))\setminus\{v_3\}\cup\{tw(v_3)\}= V(Back(P_{tw(v_3)}))$ where $P_{tw(v_3)} \in \mathcal{T}_{tw(v_3)}$.
    \end{rem}
\end{rem}
\begin{lem}\label{S_v1 = S_v2}
    For any $v_1, v_2 \in L^*(G,P,u)$, $S_{v_1} = S_{v_2}$.
    \begin{proof}
        By \cref{c = l +s}, we know that for any $v \in L^*(G,P,u)$,  $|S_v| + |L(G,P,u)| = c(v)$. If $v \in L^0(G,P,u)$, since $N(v) = N(tw(v))$, we get that $S_{v} = S_{tw(v)}.$ Thus it is enough to show that $S_{v_1} = S_{v_2}$, where $v_1, v_2 \in L(G,P,u)$. We have $S_{v_1} \cup L(G,P,u) = V(Back^*(P_{v_1})) = V(Back^*(P_{v_2})) = S_{v_2} \cup L(G,P,u)$. Therefore, we get $S_{v_1} = S_{v_2}$.
    \end{proof}
\end{lem}
\begin{notn}
    In view of \cref{S_v1 = S_v2}, we can denote set $S = S_v$, for all $v \in L^*(G,P,u)$.
\end{notn}
\vspace{1mm}
\noindent Therefore, by \cref{c = l +s}, it follows that for any $v \in L^*(G,P,u$),  we have
\begin{equation*}
    |S| = c(v) - |L(G,P,u)|
\end{equation*}
Since $Pred_w(1) \notin N(v)$, for any $v \in L^*(G,P,u)$, the pivot, $w \notin L(G,P,u)$. Thus $w \in S$ and $|S| \geq 1$.

\vspace{3mm}
\noindent What are the values that $|S|$ can take? We will analyze the two cases, $|S| =1$ and $|S| >1$. We will show that the $|S| >1$ case is not feasible in the extremal graphs.
\vspace{4mm}
\newline\textbf{$\bullet$ The structures that correspond to each  case}
\begin{lem}\label{cl:clique} (\textbf{Clique Structure})
    If $|S|=1$, then $\{w\} \cup L(G,P,u)$ induces a clique of order $c(v)$ and $L^0(G,P,u) = \emptyset$.
    \begin{proof}
        Let $v \in L(G,P,u)$, and $|S| = 1$, then $c(v) = |L(G,P,u)| +1$, thus from \cref{eq:9}, $d(v) = c(v) - 1$. On $P_v \in \mathcal{T}_v$, $v$ can only have neighbors in $Back^*(P_v)\setminus\{v\}$, otherwise we will get a cycle containing $v$ and of length more than $c(v)$. Note that, $|Back^*(P_v)\setminus\{v\}| = c(v) -1 = d(v)$. We infer that $v$ must be adjacent to all the vertices in $Back^*(P_v)\setminus\{v\}$ as shown in \Cref{clique}. Clearly, $L(G,P,u) = Back(P_v)$. As $v$ was chosen arbitrarily from $L(G,P,u)$, we get that all the vertices in $Back^*(P_v)$ are pairwise adjacent, thus $\{w\} \cup L(G, P,u)$ induces a clique. The size of the clique is $|L(G,P,u)| + 1 = c(v)$. We say that $Back^*(P_v)$ has a  \textit{Clique Structure}.
\begin{figure}[H]
    \centering
    \includegraphics[width=12cm]{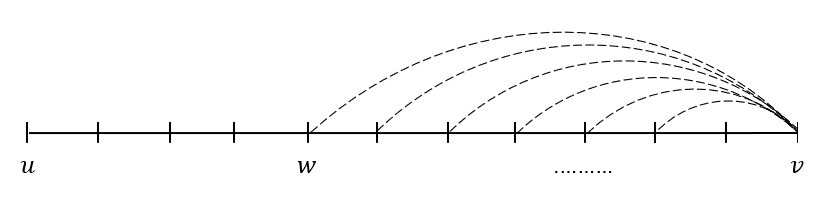}
    \caption{If $|S|=1$, then $v$ is adjacent to all the vertices of $Back^*(P_v)\setminus \{v\}$}
    \label{clique}
\end{figure}  
        \noindent Suppose there exists $\ v' \in L^0(G,P,u)$.  Then, $ tw(v') \in L(G,P,u)$. We know that $\{w\} \cup L(G,P,u)$ induces a clique. Since $tw(v')$ is adjacent to all the vertices in $\{w\}\cup L(G,P,u)\setminus\{tw(v')\}$, and $N(v') = N(tw(v'))$, $v'$ must be adjacent to $w$ and all of them. Since $|L(G,P,u)| >1$, there exists a vertex in $L(G,P,u) \setminus\{tw(v')\}$ adjacent to $v'$. This is a contradiction since a vertex in $L(G,P,u)$ can not be adjacent to any vertex outside $P_v$.

    \end{proof}
\end{lem}
\begin{lem}\label{cl:alt} (\textbf{Alternating Structure})
    If $|S|>1$, then $|S| = |L(G,P,u)|$ and every alternate vertex of $Back^*(P_v)$ starting from $w$ is a vertex in $S$, where $v \in L^*(G,P,u)$ and $P_v \in \mathcal{T}_v$.
\begin{proof}

As $|S|>1$, $S \setminus \{w\}$ is non-empty. Let $y \in S\setminus \{w\}$, be the first such vertex after $w$ on $P_v$. Since $y\in S$, we get $z_1 = Pred_y(1)$ on $P_v$ is not adjacent to $v$, since otherwise $y$ will be in $L(G,P,u)$, it follows that $z_1 \in L(G,P,u)$. On the other hand, $z_2 = Pred_y(2)\in N(v)$, on path $P_v$, since $z_2 = Pred_{z_1}(1)$ and $z_1 \in L(G,P,u)$. On path $P_v$, let $Pred_v(1) = v_1$. 
\begin{claim}\label{v_1 in S}
    $v_1 \in S$
    \begin{proof}
        Consider path $P_{z_1} \in \mathcal{T}_{z_1}$, as in \Cref{alt1}, formed by removing edge $(z_2,z_1)$ and adding the edge $(z_2,v)$ to $P_v$. We know that $P_{z_1}$ does not have holes. Since $v\notin N(z_1)$, arguing with respect to $P_{z_1}$ we infer that, $v_1 \notin L(G,P,u)$. Therefore we get $v_1 \in S$.
       \end{proof}
\end{claim}
\begin{figure}[H]
    \centering
    \includegraphics[width=12cm]{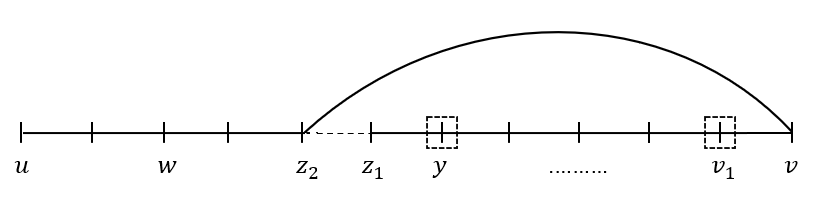}
    \caption{Boxed vertices represent the vertices in $S$. For $|S|>1$, $y$ is assumed to be the first vertex in $S\setminus \{w\}$ on $P_v$. With respect to $P_v$, $z_1 = Pred_y(1)$, $z_2 = Pred_y(2)$ and $v_1 = Pred_v(1)$.}
    \label{alt1}
\end{figure}  
\begin{claim}\label{z_2 = w}
    $z_2 = w$ 
    \begin{proof}
        As $w \in N(v)$, we have $Succ_w(1) \in L(G,P,u)$. But $y\in S\setminus\{w\}$ and thus, $y\neq Succ_w(1) $. We infer that $z_1 \neq w$. Suppose $z_2 \neq w$. Then there exists a vertex $z_3 = Pred_y(3)$ in $Back^*(P_v)$ before $z_2$. As $y$ is the first element of $S\setminus \{w\}$ on $P_v$, $z_1, z_2 \in L(G,P,u)$ and therefore $z_2,z_3 \in N(v)$.\newline Now consider the path $P_{z_2}\in \mathcal{T}_{z_2}$, as shown in \Cref{alt2}, formed by removing the edge $(z_3,z_2)$ and adding the edge $(z_3,v)$ to $P_v$. Since $v \in N(z_2)$, arguing with respect to path $P_{z_2}$, we have $v_1 \in L(G,P,u)$. Thus $v_1 \notin S$. This contradicts \cref{v_1 in S}. Therefore, $z_2 = w$.
    \end{proof}
    \begin{figure}[H]
    \centering
    \includegraphics[width=12cm]{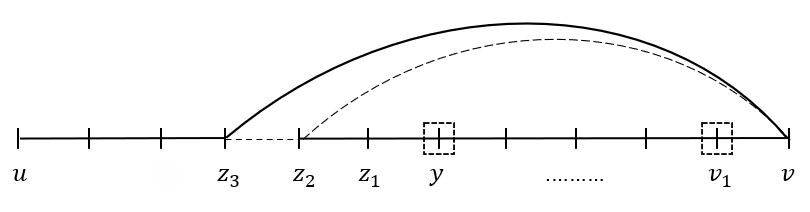}
    \caption{Path $P_{z_2}$ can be constructed if $z_3$ exists in $Back^*(P_v)$ before $z_2$.}
    \label{alt2}
\end{figure} 
    
\end{claim}
\noindent Since $y \in N(v)$, $y_1 = Succ_y(1) \in L(G,P,u)$ on path $P_v$. Consider the path $P_{y_1}$, formed by deleting the edge $(y,y_1)$ and adding edge $(y,v)$ to the path $P_v$. On the path $ P_{y_1}$ since $v_1 \in S$, it follows that $v \notin N(y_1)$; otherwise, we would have $v_1 \in L(G, P, u)$ contradicting \cref{v_1 in S}.
 Now considering path $P_v$, since we have $y_1\notin N(v)$, we get $Succ_y(2) \in S$. Following a sequence of similar arguments, one can show $Succ_y(4), Succ_y(6), \dots $ and so on are in $S$. Note that since $v_1 \in S$, we get $Pred_{v_1}(1) \notin N(v).$ Therefore, $Pred_{v_1}(1) \notin S$.
\newline Thus, $\{w = Pred_y(2), y , Succ_y(2) , Succ_y(4), \dots ,v_1 = Pred_v(1)\} = S$. Since on $P_v$, each vertex in $S$ is followed by a vertex in $L^*(G,P,u)$, and from \cref{|L| = |L_v|} we have $|Back^*(P_v)\cap L^*(G,P,u)| = |L(G,P,u)|$, we get $|S| = |L(G,P,u)|$. Thus, the set $S$, as depicted in \Cref{alt}, consists precisely of the vertices at even distances from the pivot $w$ within $Back^*(P_v)$.

\end{proof}
 \begin{figure}[H]
    \centering
    \includegraphics[width=12cm]{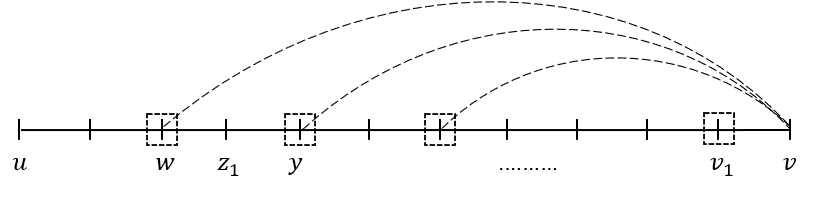}
    \caption{Every alternate vertex of $Back^*(P_v)$ starting from $w$ is a vertex in $S$.}
    \label{alt}
\end{figure} 
\end{lem} 

\begin{lem}\label{S=N(x)}
    If $P_v$ exhibits alternating structure, then $S = N(x)$ for all $x \in L^*(G,P,u)$.
    \begin{proof}
         Recall that from \cref{S_v1 = S_v2}, $S = S_x \subseteq N(x)$ for all $x \in L^*(G,P,u)$. From \cref{eq:9} we have, for all $x \in L^*(G,P,u)$, $d(x) = |L(G,P,u)|$. Since $P_v$ exhibits alternating structure, from \cref{cl:alt} we have $|S| = |L(G,P,u)| = d(x) = |N(x)|$. Therefore, $S = N(x)$ for all $x \in L^*(G,P,u)$.
    \end{proof}
\end{lem}
\noindent We will show that the \textit{alternating structure} as described in \cref{cl:alt}, is not feasible if there is equality in the bound of \Cref{th:main}. Until otherwise stated, we assume $|S| >1$ while exploring the alternating structure of the path $P$.
\vspace{4mm}
\newline\textbf{$\bullet$ Alternating Structure is not feasible}

\begin{figure}[H]
    \centering
    \includegraphics[width=12cm]{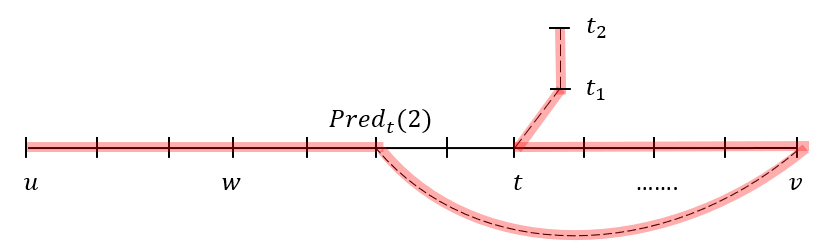}
    \caption{The highlighted path, $P'$, is a longer $u$-path than $P_v$}
    \label{2_branch}
\end{figure}

\begin{figure}[H]
    \centering
    \includegraphics[width=12cm]{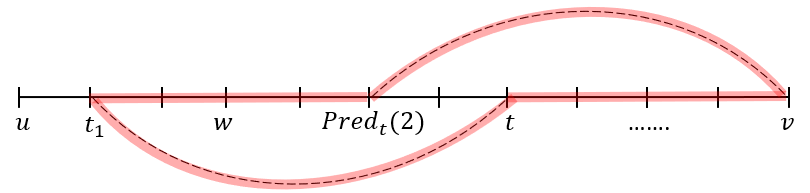}
    \caption{The length of the highlighted cycle $Q$, is at least $c(v) +1$ and it contains $v$ }
    \label{Adjacent}
\end{figure}

\begin{notn}
    For a vertex $x \in V(P_v)$, the path $Q = v, t_1, t_2$ is referred to as a \textit{$2$-branch} attached to $x$ if $t_1, t_2 \notin V(P)$.

\end{notn}
\begin{lem}\label{cl:branch}
    No vertex from $S\setminus \{w\}$ can have a $2$-branch attached to it.
    \begin{proof}
        Suppose there exists $t \in S\setminus\{w\}$ such that there exists a $2$-branch, $Q = t,t_1,t_2$. As $t \neq w$, $Pred_t(2)\in S$. Thus, $Pred_t(2) \in N(v)$. Now construct a path $P'$ from path $P_v$, by removing the vertex $Pred_t(1)$ and adding the edges $(Pred_t(2), v), (t,t_1), (t_1,t_2)$, that is \[P' = P_v(u,Pred_t(2)),v,P_v(v,t),t_1,t_2\] Refer to the highlighted path in \Cref{2_branch}. \newline Thus $P'$ is a $u$-path with $|V(P)| +1$ vertices. This contradicts the fact that $P$ is a longest $u$-path.
    \end{proof}
\end{lem}
\begin{lem}\label{cl:nbr}
    No vertex from $S\setminus \{w\}$ can be adjacent to any vertex on $P_v(u,Pred_w(2))$.
    \begin{proof}
        Suppose there exists $t \in S\setminus\{w\}$ adjacent to some vertex, $t_1 \in V(P_v(u,Pred_w(2)))$. On path $P_v$, since $t \neq w$ we get $Pred_t(2) \in S$. Thus, $Pred_t(2) \in N(v)$. Construct a cycle,\[Q = P_v(t_1,Pred_t(2)),v,P_v(v,t),t_1\] by adding the edges $(t_1,t)$ and $(Pred_t(2),v)$, as highlighted in \Cref{Adjacent}.\newline $Q$ contains all the vertices on $P_v$ from $Pred_w(2)$ to $v$ except $Pred_t(1)$. Thus $|V(Q)| \geq c(v) + 1$. So, $Q$ is a cycle containing $v$ and of size at least $c(v) +1$ resulting in a contradiction.
    \end{proof}
\end{lem}

\noindent
For equality in \Cref{th:main}, from \cref{c_i(v) = c(v)}, we infer that the removal of $L^*(G,P,u)$ must not lead to a reduction in the weights of the vertices in the graph $G_1 = G \setminus L^*(G,P,u)$, except possibly for $u$. In graph $G$ vertices of $S \setminus\{w\}$ are part of a cycle of length $c(v)$. Therefore the vertices in $S \setminus \{w\}$ must still belong to a cycle of length at least $c(v)$ in $G_1$.
\newline Let $C$ denote the longest cycle in $G_1$ that includes at least one vertex from $S \setminus \{w\}$. We now proceed to analyze the structure of $C$.

\vspace{1mm}
\noindent\textbf{Case 1:}\label{case 1} There exists at least two consecutive vertices on $C$ that are not contained in $S$.

\vspace{2mm}
\noindent Since $C$ contains at least one vertex from $S\setminus\{w\}$ and at least two consecutive vertices that are not from $S$,  it is easy to see that we can find a vertex $t \in S\setminus\{w\}$ and $p_1,p_2\notin S$ such that $t, p_1,p_2$ are consecutive on $C$. Since $t\neq w$, with respect to path $P_v$ we get, $Pred_t(2) \in S$.

\vspace{2mm}
\noindent\textbf{Subcase 1.1.} $p_1 \notin V(P_v)$.

\begin{figure}[H]
    \centering
    \includegraphics[width=12cm]{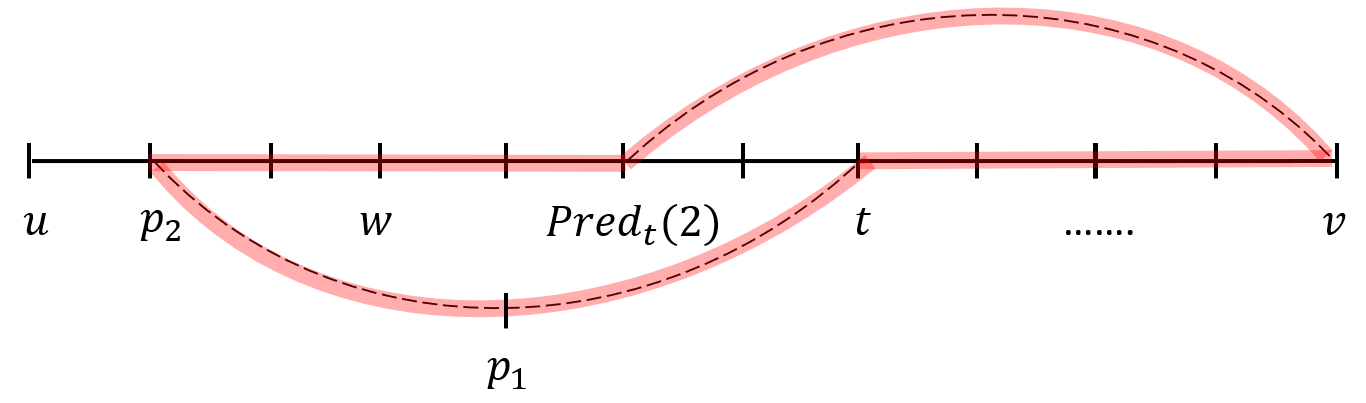}
    \caption{The length of the highlighted cycle $Q$, is at least $c(v) +1$ and it contains $v$ }
    \label{sbc1.1}
\end{figure}

\vspace{2mm}
\noindent Since $p_2 \in V(G_1)$, it follows that $p_2 \notin L^*(G, P, u)$. By the assumption case 1, $p_2 \notin S$. From \Cref{cl:branch}, we conclude that $p_2 \in V(P)$; otherwise, a $2$-branch would be attached to $t$. But $p_2 \notin L^*(G,P,u)\cup S = Back^*(P_v)$. Therefore, $p_2 \in V(Front(P_v))$. Recall that since $t \in S\setminus \{w\}$, $Pred_t(2) \in S$ on path $P_v$. The highlighted cycle, \[Q=P_v(p_2,Pred_t(2)),v,P_v(v,t),p_1,p_2\] in \Cref{sbc1.1}, contains all the vertices from $P_v(Pred_w(1),v)$ except the vertex $Pred_t(1)$. Note that $Q$ contains $p_1\notin V(P_v)$ also. Thus, $|V(Q)|\geq c(v) +1$. So, $Q$ is a cycle containing $v$ and of size at least $c(v) +1$ resulting in a contradiction.

\vspace{2mm}
\noindent\textbf{Subcase 1.2.} $p_1 \in V(P_v)$.

\vspace{2mm}
\noindent By the assumption of case 1, we have $p_1 \notin S$. Using \Cref{cl:nbr}, we can say that $p_1 \notin V(P_v(u,Pred_w(2)))$. Since $p_1\in V(G_1)$, it follows that  $p_1 \notin L^*(G,P,u)$. Therefore, $p_1 \notin L^*(G,P,u)\cup S\cup V(P_v(u,Pred_w(2))) = V(P_v) \setminus \{Pred_w(1)\}$. Thus $p_1 = Pred_w(1)$.
\vspace{2mm}
\newline \textbf{Claim:} $C$ must contain a vertex from $S\setminus\{t\}$. 
\begin{proof}
    Suppose not. Let $t,x,y$ be consecutive vertices on $C$ where $x \neq p_1$. By assumption, $x \notin S$ and $x \notin L^*(G,P,u)$. From \cref{cl:nbr}, we get that $x \notin V(P_v(u,Pred_w(2)))$. Recall that $p_1 = Pred_w(1)$ and $x \neq p$, from this we infer that, $x \notin V(P_v)$. Also note that, $y \notin S\cup L^*(G,P,u)$. Since $|S| >1$, we have $|C| \geq 2|S| \geq 4$. Therefore, $y \neq p_1$. From \cref{cl:branch}, we get that $y \in V(P_v)$; otherwise we get a 2-branch attached on $t$. Thus $y \in V(P_v(u,Pred_w(2)))$. Now consider the cycle $Q_1$, highlighted in \cref{claim1.2.1};
    \[Q_1 = P_v(y,Pred_t(2)),v,P_v(v,t),x,y\]
    \begin{figure}[H]
    \centering
    \includegraphics[width=12.8cm]{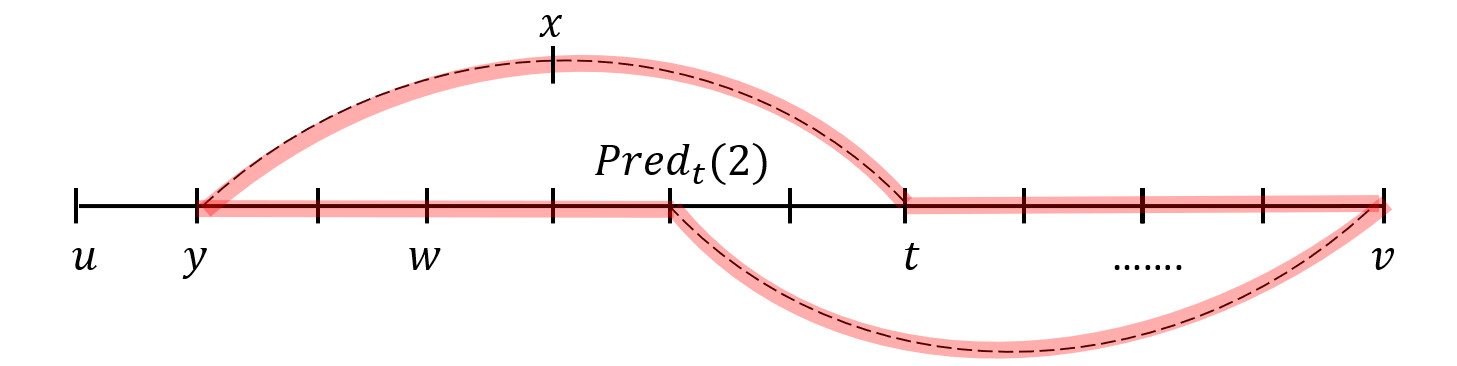}
    \caption{The length of the highlighted cycle $Q_1$, is at least $c(v) +1$ and it contains $v$ }
    \label{claim1.2.1}
\end{figure} 
    \noindent $Q_1$ contains all the vertices from $P_v(y,v)$ except the vertex $Pred_t(1)$. Note that $Q_1$ contains $x \notin V(P_v)$ also. Thus, $|V(Q_1)| \geq c(v)+2$. So, $Q_1$ is a cycle containing $v$ and of size at least $c(v)+2$ resulting in a contradiction.
\end{proof}
\noindent\textbf{Claim: }By the case 1 assumption, $p_2 \notin S$. Let $C = t,p_1,p_2,v_1,v_2,\dots, v_k,t$. Suppose $ j = min\{\ i\in [k]:\ v_i \in S\}$, which exists due to the above claim. Then $v_j = w$.
\begin{proof}
    Suppose $v_j \neq w$. Clearly, on path $P_v$, $Pred_{v_j}(2) \in S$. Recall that $p_1 = Pred_w(1)$, for subcase 1.2. Moreover $p_2,v_1,\dots,v_{j-1} \notin Back^*(P_v)$, since $V(C) \cap L^*(G,P,u) = \emptyset$ and $v_j$ is the first vertex from set $S$ after $t$ on $C$ in this direction. Moreover $P_v(v_j,v) \subseteq Back^*(P_v)$ and since $Pred_{v_j}(2) \in S$ we have, $p_2,v_1,\dots,v_{j-1} \notin P_v(Pred_{v_j}(2),p_1)$. Therefore we can define cycle $Q_2$ as shown in \cref{1.2claim3}
    \[Q_2 = p_1,p_2,v_1,\dots,v_j,P_v(v_j,v),Pred_{v_j}(2),P_v(Pred_{v_j}(2),p_1)\]
    \begin{figure}[H]
    \centering
    \includegraphics[width=12cm]{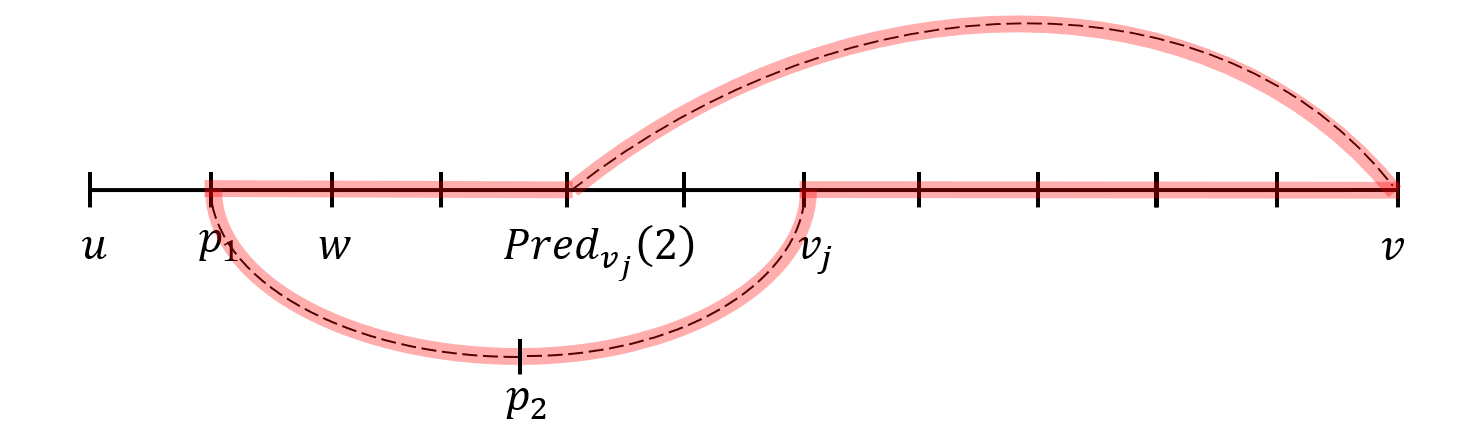}
    \caption{The highlighted cycle $Q_2$, has length at least $c(v) +1$ and contains vertex $v$. The vertex $p_2$ may appear on $P_v$ before $p_1$.}
    \label{1.2claim3}
\end{figure} 
    \noindent $Q_2$ contains all the vertices from $P_v(p_1, v)$ except the vertex $Pred_{v_j}(1)$ and also contains the vertex $p_2 \notin P_v(p_1,v)$. Thus, $|V(Q_2)| \geq c(v) +1$. So, $Q_2$ is a cycle containing $v$ and of size at least $c(v) +1$ resulting in a contradiction.
\end{proof}

\begin{figure}[H]
    \centering
    \includegraphics[width=12cm]{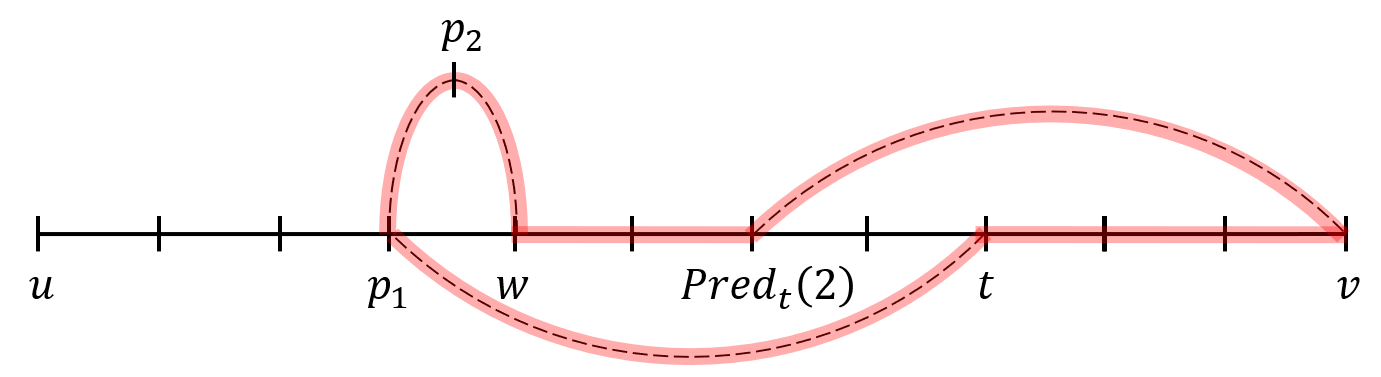}
    \caption{The highlighted cycle $Q_3$, has length at least $c(v) +1$ and contains vertex $v$. The vertex $p_2$ may appear on $P_v$ before $p_1$.}
    \label{sbc1.2}
\end{figure}

\noindent We have $v_j = w$. Consider the cycle $Q_3$, as shown in the \cref{sbc1.2};
\[Q_3 = p_1,p_2,v_1,\dots,v_{j-1},w,P_v(w,Pred_t(2)),v,P_v(v,t),p_1\]
$Q_3$ contains all the vertices from $P_v(p_1,v)$ except the vertex $Pred_t(1)$ and also contains the vertex $p_2 \notin P_v(p_1,v)$. Thus $|V(Q_3)| \geq c(v) + 1$. So, $Q_3$ is a cycle containing $v$ and of size at least $c(v) +1$ resulting in a contradiction.

\noindent\textbf{Case 2:}\label{case 2}
    For every pair of consecutive vertices on $C$, at least one from the pair is in $S$.
    
\begin{figure}[H]
\centering
    \includegraphics[width=4.2cm]{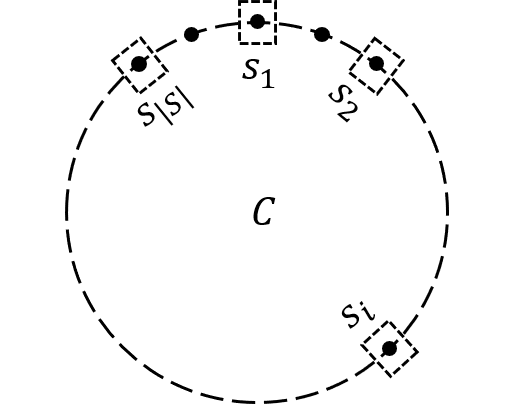}
    \caption{$C$ is a $2|S|$ length cycle. The vertices of $S$ (boxed) are positioned alternately on $C$.}   \label{case2}

\end{figure}

    \vspace{2mm}
    \noindent From \cref{c_i(v) = c(v)} we get that the cycle $C$ has a length of at least $c(v) = 2|S|$, that is $|V(C)| \geq 2|S|$. From the assumption of Case 2, clearly $|V(C)| \leq 2|S|$. Therefore, $|V(C)| = 2|S|$. From this, we can conclude that, $V(C)$ must include all the vertices from $S = \{s_1,s_2, \dots, s_{|S|}\}$, arranged alternately along $C$, as represented by the boxed vertices in \cref{case2}.

\noindent\textbf{Subcase 2.1.} $Pred_w(1)\in V(C)$

\begin{figure}[H]
\centering
    \includegraphics[width=12.2cm]{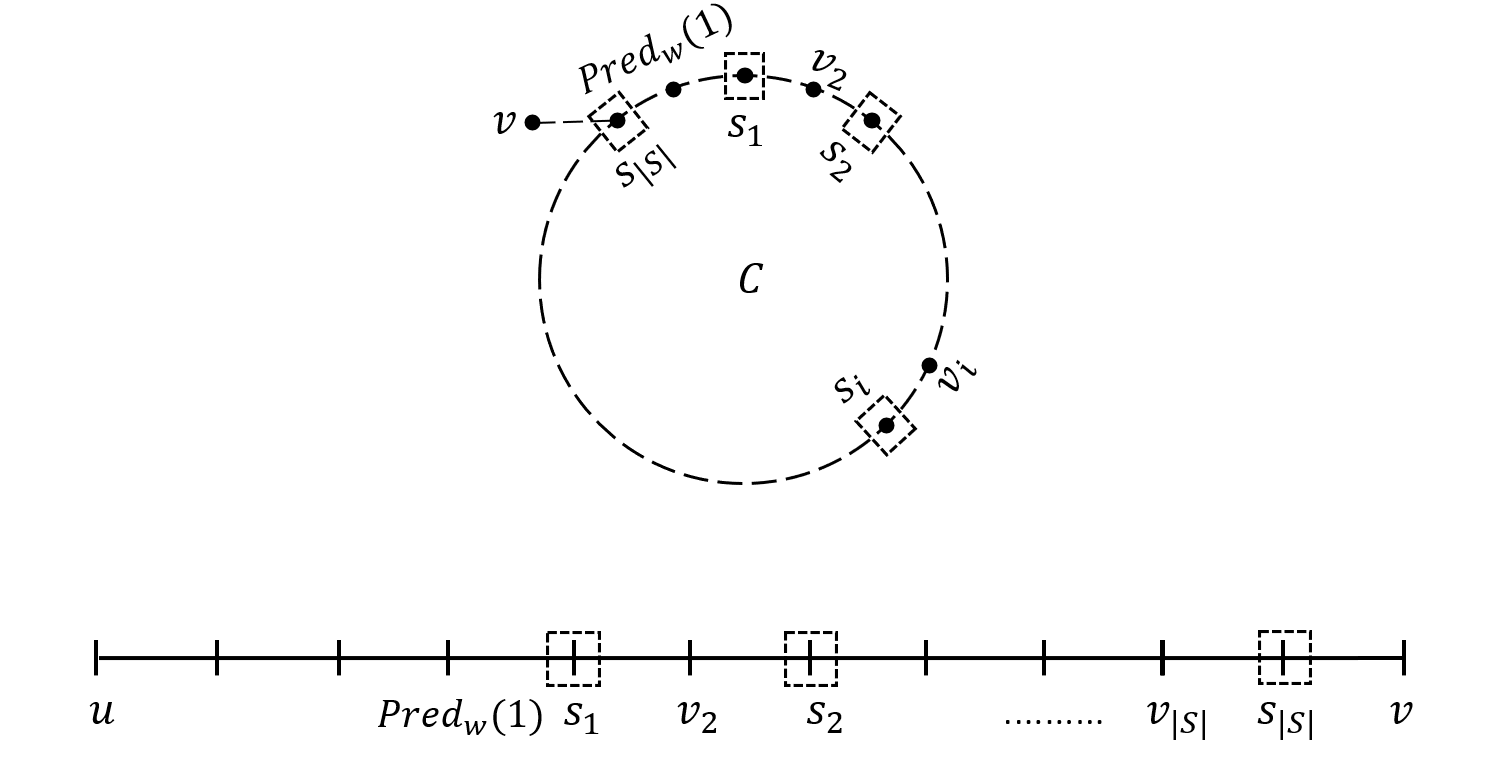}
    \caption{$C$ contains $Pred_w(1)$ and the path $P'$ is created in the cyclic order}
    \label{sbc2.1}

\end{figure}

\vspace{2mm}
\noindent Without loss of generality let $C = Pred_w(1),s_1,v_2,s_2,v_3,\dots, v_{|S|}, s_{|S|}, Pred_w(1)$. Here $S = \{s_1,s_2, \dots, s_{|S|}\}$. Clearly $v_i \notin Back^*(P_v)$. From \cref{cl:nbr} we get that $v_i \notin V(P_v(u,Pred_w(2)))$, and $v_i \neq Pred_w(1)$ for all $i$. Therefore, $v_i \notin V(Back^*(P_v))\cup V(P_v(u,Pred_w(2))) \cup \{Pred_w(1)\} = V(P_v)$. 
\newline  Recall that $S \subseteq  N(v)$. Therefore, $s_{|S|}$ is adjacent to $v$ in $G$. Define the $u$-path $P'$, in graph $G$, as shown in \Cref{sbc2.1}
\[P' = Front(P_v),s_1,v_2,s_2,\dots,v_{|S|},s_{|S|},v\]
Both $P_v$ and $P'$ have the same length; therefore, $P'$ could have been chosen in place of $P_v$ from the outset. In this case, for equality to hold in \Cref{th:main}, $P'$ must exhibit either the \textit{clique structure} as in \cref{cl:clique} or the \textit{alternating structure} as specified in \cref{cl:alt}.
\noindent From \cref{S=N(x)} we have $N(v) = S$. Thus, $v$ has no neighbor outside path $P'$, and on $P'$ the farthest neighbor of $v$ in distance is $s_1$. Considering path $P'$, clearly $dist(s_1,v) = c(v) -1$. Thus, $s_1$ is the pivot of the path $P'$. \vspace{2mm}\newline If $P'$ exhibits the clique structure, the corresponding set to be removed, denoted as $L^*(G,P',u)$, would be $\{v_2, s_2,v_3, \dots,v_{|S|}, s_{|S|}, v\}$. It is easy to see that $L^*(G,P',u)\cap L^*(G,P,u) = \{v\}$, since $L^*(G,P',u)\setminus \{v\} \subseteq V(C) \subseteq V(G)\setminus L^*(G,P,u)$. In \cref{|L| >1}, we assumed that $|L(G,P,u)| \geq 2$, therefore $L^*(G,P,u) \setminus \{v\} \neq \emptyset$. Also, we assumed $|S| \geq 2$, therefore $s_2 \in L^*(G,P',u)$.  Let $x \in L^*(G,P,u)\setminus \{v\}$, then $x \notin V(P')$. But $s_2$ is adjacent to $x$. This results a contradiction as $L^*(G,P',u)$ can not have a neighbor outside $P'$; otherwise, one can get a $u$-path longer than $P'$ and thus longer than $P_v$.

\vspace{2mm}
\noindent Thus, $P'$ must have alternating structure. Therefore $L^*(G,P',u) = \{v_2,\dots,v_{|S|},v\}$. Using 
 \cref{S=N(x)}, we have $N(v_i) = S$ for all $i$. Since $P_v$ exhibits alternating structure, we already have $N(x) = S$ for all $x \in L^*(G,P,u).$  Since $|S| \geq 2$, $\{v_2,\dots,v_{|S|}\} \neq \emptyset$. Note that, for all $y \in \{v_2,v_3,\dots,v_{|S|}\}$ there exists $x \in L(G,P,u)$ such that $N(y) = N(x) (=S)$. Therefore, $\{v_2,\dots,v_{|S|}\} \subseteq L^0(G,P,u) \subseteq L^*(G,P,u)$. Note that $v_i \in V(G_1)$ thus $v_i \notin L^*(G,P,u)$. Thus we arrive at a contradiction.

\vspace{2mm}
\noindent\textbf{Subcase 2.2.} $Pred_w(1) \notin V(C)$

\vspace{2mm}
\noindent Let $C = s_1, v_1, s_2, v_2, \dots, s_{|S|}, v_{|S|}, s_1$, where $S = \{s_1, s_2, \dots, s_{|S|}\}$. From \Cref{cl:nbr}, all the $v_i$'s lie outside $P_v$. Without loss of generality, assume $s_1 = w$. Then, define the $u$-path  
\[
P' = \text{Front}(P_v),s_1, v_1, s_2, v_2, \dots, s_{|S|}, v_{|S|}.
\]  
Since $P'$ is a $u$-path of the same length as $P_v$, it follows that $P'$ could have been chosen initially instead. Suppose $v_{|S|}$ is adjacent to some vertex,  $t \in Front(P_v)$. Then we can get a cycle $C'$, containing the vertices of $C$ and those between $t$ and $s_1$ on $P'$, including $t$. Note that $Front(P_v) \subseteq V(G_1)$, thus $V(C') \subseteq V(G_1)$. This contradicts the fact that $C$ is the longest cycle in $G_1$ which includes at least one vertex from $S\setminus \{w\}$. Therefore, $v_{|S|}$ is not adjacent to any vertex in $Front(P_v)$. For equality in \Cref{th:main}, $P'$ either exhibits the clique structure as in \cref{cl:clique} or possesses the alternating structure as shown in \cref{cl:alt}. Clearly, $v_{|S|}$ is not adjacent to any vertex outside $P'$; otherwise, we get a longer $u$-path than $P_v$. Also $v_{|S|}$ does not have any neighbor in $Front(P_v)$. 

\vspace{2mm}
\noindent Suppose $P'$ exhibits clique structure then $L^*(G,P',u) = V(P')\setminus Front^*(P_v) = \{v_1,s_2,\dots,s_{|S|},v_{|S|}\}$, since $s_1$ is the farthest neighbor of $v_{|S|}$ on $P'$. Applying similar reasoning as subcase 2.1, we get $s_2$ to be adjacent to $v \notin V(P')$. Thus resulting in a longer $u$-path than $P'$, therefore contradicting that $P'$ (or $P_v$) is the longest $u$-path in $G$.

\vspace{2mm}
\noindent Thus, $P'$ must have alternating structure. Since $s_1$ is the first neighbor of $v_{|S|}$ on $P'$, $L^*(G,P',u) = \{v_1,v_2,\dots,v_{|S|}\}$. From \cref{S=N(x)}, we have $N(v_i) = S$ for all $i \in [|S|]$. But $P_v$ exhibits alternating structure, therefore $N(x) = S$ for all $x \in L^*(G,P,u)$. Clearly, $v_i \notin V(P_v)$, and for any $v_i$ there exists $x \in L(G,P,u)$ such that $N(x) = N(v_i) (= S)$. Therefore, $v_i \in L^0(G,P,u) \subseteq L^*(G,P,u)$. This contradicts the fact that $v_i \in V(G_1)$, and thus $v_i\notin L^*(G,P,u)$.

\vspace{2mm}
\noindent We deduce that $P'$ cannot possess either of the configurations described in \cref{cl:clique} or \cref{cl:alt}, namely the clique structure and the alternating structure, respectively.

\vspace{2mm} \noindent Consequently, neither Case 1 nor Case 2 can occur under equality. This implies that the alternating structure described in \cref{cl:alt} is not feasible. Therefore, for equality, the longest $u$-path in $G$ must exhibit the clique structure, as outlined in \cref{cl:clique}.
\vspace{2mm}
\newline\textbf{$\bullet$ Clique Structure implies that  $G$ is a Parent-dominated block graph}
\begin{lem}
    If $G$ is an extremal graph for \Cref{th:main}, then $G$ is a parent-dominated block graph.
    \begin{proof}
        Proof by Induction on the number of vertices:

\noindent\textit{Base Case:} If $|V(G)| = 1$, then $G$ is trivially a parent-dominated block graph.

\noindent\textit{Induction Hypothesis:} Assume that any extremal graph for \Cref{th:main} with fewer than $|V(G)|$ vertices is a parent-dominated block graph.

\noindent\textit{Induction Step:} Let $G$ be an extremal graph. If $G$ is a clique, we are done. Thus, assume $G$ is not a clique. By the induction hypothesis, the graph $G_1 = G \setminus L^*(G, P, u)$ is a parent-dominated block graph. $G$ is connected by \cref{connected}. We know that $Back^*(P) = L^*(G, P, u) \cup \{w\}$ induces a clique in $G$ by \cref{cl:clique}, where $w$ is the only vertex in $G_1$ that is adjacent to the vertices in $L^*(G, P, u)$. Thus, $G$ is a block graph, with $B = Back^*(P)$ as a leaf block and ${w}$ as its cut vertex. Let $B'$ be the parent block of $B$ in $G$. Thus $V(B) \cap V(B') = \{w\}$. The weight of $w$ should remain unchanged upon the removal of $L^*(G, P, u)$ by \cref{c_i(v) = c(v)}. Since $G_1$ is a parent-dominated block graph, $c(w) = |V(B')|$. As the vertices in $B$ induce a clique and $w \in V(B)$, we get $c(w) \geq |V(B)|$. It follows that $|V(B')| \geq |V(B)|$. This implies that $G$ is a parent-dominated block graph.
\vspace{2mm}
\newline Therefore, if $G$ is an extremal graph for \Cref{th:main} with $\delta(G) \geq 2$, then $G$ is a parent-dominated block graph. Now, suppose $\delta(G) = 1$ and $G$ is extremal. We can iteratively remove degree-one vertices until we obtain a graph $H$ with $\delta(H) \geq 2$. By \cref{|L| >1}, $H$ remains extremal for \Cref{th:main}; it must also be a parent-dominated block graph. The removed degree-one vertices act as leaf blocks (at the times they are removed), ensuring that $G$ is still a block graph. Since these degree-one vertices form blocks of size $2$ (edges), the parent-domination property continues to hold. Thus, we conclude that $G$ is a parent-dominated block graph.

    \end{proof}
\end{lem}

\subsection{Proof of \Cref{th:main:path}}\label{Proof of 2.1}
\begin{proof}We will prove \Cref{th:main:path} using induction on the number of vertices in graph $G$.

    \vspace{1mm}
    \noindent \textit{Base Case:}  
    If $|V(G)| =1$, let $V(G) = \{v\}$. Clearly $|E(G)| = 0 = \frac{p(v)}{2}$. $G \cong K_1$ and equality holds.

    \vspace{1mm}
    \noindent \textit{Induction Hypothesis:} Assume that for any graph with at most $|V(G)| -1$ vertices, the inequality in \Cref{th:main:path} is true. And equality holds if and only if all the components of the graph are cliques.

    \vspace{1mm}
    \noindent \textit{Induction Step:} If $G$ is disconnected by induction hypothesis the statement of the theorem is true for each connected component and thus for $G$. Therefore, we assume that $G$ is connected. Let $k = max\{ p(v)\mid v\in V(G)\}$. The discussion will proceed by considering separate cases.

\vspace{2mm}
\noindent\textbf{Case 1:}
     There exists a path $P$ in $G$ of length $k$ such that the endpoints of $P$ are adjacent. 
 \begin{claim}
     $V(P) = V(G)$.
     \begin{proof}
         Let $P = v_0,v_1,\dots,v_k$. By assumption of the case, we get that $v_0$ is adjacent to $v_k$. Suppose $V(G) \setminus V(P) \neq \emptyset$. Since $G$ is connected, there exists a vertex  $v \in V(G)\setminus V(P)$, adjacent to some vertex in $P$. Since the vertices of $P$ form a cycle, without loss of generality, we can assume $v$ is adjacent to $v_0$. Now consider the path;
         \begin{equation*}
             Q = v,v_0,v_1,\dots,v_{k}
         \end{equation*}
Thus $Q$ is a path in $G$ with length $k +1$. Thus we arrive at a contradiction to the assumption that the length of a longest path in $G$ is $k$.
     \end{proof}
 \end{claim}
\noindent Since $V(P) = V(G)$, we get that $p(v) = |V(G)| -1 = (n-1)$, for all $v \in V(G)$. The maximum number of edges $G$ can have is $\frac{n(n -1)}{2}$. Thus we get that;
\begin{equation}\label{case1path}
    |E(G)| \leq \frac{n(n-1)}{2} = \sum_{v\in V(G)}\frac{p(v)}{2}
\end{equation}
Clearly, equality in \cref{case1path} holds if and only if $G$ is a clique.
\vspace{2mm}
\newline\textbf{Case 2:}
    If $P$ is a path of length $k$ in $G$, then the endpoints of $P$ are not adjacent. That is there are no cycles in $G$ with length $k+1$.
\vspace{1mm}    

\noindent Let $P = v_0,v_1,\dots,v_k$.  Let $d(x)$ be the degree of the vertex $x \in V(G)$. Without loss of generality assume $d(v_0) \geq d(v_k)$.
\begin{claim}
    $d(v_k) \leq \frac{k}{2}$
    \begin{proof}
        Suppose not, that is $d(v_k) \geq \frac{k+1}{2}$. Thus we get $d(v_0) \geq \frac{k+1}{2}$. Since $P$ is a longest path in $G$ and $v_0$ is not adjacent to $v_{k}$, we get, $N(v_0), N(v_k) \subseteq \{v_1,v_2,\dots, v_{k-1}\}$. Let $N = N(v_k) \setminus \{v_{k-1}\}$ and $N^+ = \{v_{j+1} \mid v_j \in N\}$. Clearly, $|N| = |N^+| \geq \frac{k+1}{2}-1 = \frac{k-1}{2}$. Note that $v_0, v_k \notin N^+$. Therefore, $|\{v_1,v_2,\dots,v_{k-1}\} \setminus N^+| \leq \frac{k-1}{2}$. Since $|N(v_0)| \geq \frac{k+1}{2}$, by pigeonhole principle, we get that, there exists $v_i \in N^+$ such that $v_i \in N(v_0)$ and $v_{i-1} \in N(v_k)$. Consider the cycle; $$K = v_0, v_1, \dots, v_{i-1},v_k,v_{k-1},\dots, v_i,v_0$$ Clearly $K$ is a $k+1$ length cycle in $G$. Thus we get a contradiction.
    \end{proof}
    \noindent Let $G_1 = G \setminus \{v_k\}$ and $p_{G_1}(v)$ denote the weight of a vertex in $G_1$. Clearly $p(v) \geq p_{G_1}(v)$ for all $v \in V(G_1)$. Using the induction hypothesis we get;
    \begin{equation}\label{path_induction}
        |E(G_1)| \leq \sum_{v \in V(G_1)}\frac{p_{G_1}(v)}{2} \leq \sum_{v \in V(G_1)}\frac{p(v)}{2}
    \end{equation}
    Note that, $|E(G)| = |E(G_1)| + d(v_k)$. Thus from \cref{path_induction} we get;
    \begin{equation}\label{path_ineq}
        |E(G)| \leq \sum_{v \in V(G_1)}\frac{p(v)}{2} + d(v_k) \leq \sum_{v \in V(G_1)}\frac{p(v)}{2} + \frac{p(v_k)}{2} = \sum_{v \in V(G)}\frac{p(v)}{2}
    \end{equation}
\end{claim}

\noindent For equality in \cref{path_ineq}, we must have equality in \cref{path_induction}. By induction hypothesis, the first inequality of \cref{path_induction} becomes an equality if and only if all the components of $G_1$ are cliques. Since the vertices $\{v_0,v_1,\dots,v_{k-1}\}$ form a path in $G_1$, they belong to the same connected component of $G_1$, say $C$. For equality in the second inequality of \cref{path_induction} we need to have $p_{G_1}(v) = p(v)$ for all $v \in V(G_1)$. Thus, $p_{G_1}(v_i) = k$ for all $i \in \{0,1\dots, k-1\}.$ Since $C$ is a clique with vertices having weight $k$, $C$ has order $k+1$, that is, there is one more vertex in $C$ other than $\{v_0,\dots,v_{k-1}\}$. Thus, $C$ contains a cycle of length $k+1$. This contradicts the assumption that $G$ has no $k+1$ length cycle.
\vspace{2mm}
\newline Thus equality is not feasible in Case 2. Therefore, we have equality in \Cref{th:main:path} if and only if all the components of $G$ are cliques.

\end{proof}

\bibliographystyle{plain}
\bibliography{references}
\end{document}